\documentclass[11pt,a4paper]{article}

\usepackage{amsmath,amssymb,latexsym,pstricks,mathrsfs,comment,amsthm,graphicx,tikz,enumerate,
pgffor,cite,wrapfig,float,url}
\usepackage{bibspacing}
\newcommand{\nc}{\newcommand}
\nc{\rnc}{\renewcommand}

\usepackage{a4wide}  

\rnc{\ss}{\smallskip} \nc{\ms}{\medskip} \nc{\bs}{\bigskip}
\nc{\bsa}{\vspace{-\bigskipamount}}

\makeatletter

\begin{document}

\hyphenation{mon-oid mon-oids}

\nc{\Sone}{S^1}
\nc{\keru}{\operatorname{ker}^\wedge} \nc{\kerl}{\operatorname{ker}_\vee}

\nc{\coker}{\operatorname{coker}}
\nc{\KER}{\ker}
\nc{\N}{\mathbb N}
\nc{\LaBn}{L_\al(\B_n)}
\nc{\RaBn}{R_\al(\B_n)}
\nc{\LaPBn}{L_\al(\PB_n)}
\nc{\RaPBn}{R_\al(\PB_n)}
\nc{\rhorBn}{\rho_r(\B_n)}
\nc{\DrBn}{D_r(\B_n)}
\nc{\DrPn}{D_r(\P_n)}
\nc{\DrPBn}{D_r(\PB_n)}
\nc{\DrKn}{D_r(\K_n)}
\nc{\alb}{\al_{\vee}}
\nc{\beb}{\be^{\wedge}}
\nc{\bnf}{\bn^\flat}
\nc{\Bal}{\operatorname{Bal}}
\nc{\Red}{\operatorname{Red}}
\nc{\Pnxi}{\P_n^\xi}
\nc{\Bnxi}{\B_n^\xi}
\nc{\PBnxi}{\PB_n^\xi}
\nc{\Knxi}{\K_n^\xi}
\nc{\C}{\mathbb C}
\nc{\exi}{e^\xi}
\nc{\Exi}{E^\xi}
\nc{\eximu}{e^\xi_\mu}
\nc{\Eximu}{E^\xi_\mu}
\nc{\REF}{ {\red [Ref?]} }
\nc{\GL}{\operatorname{GL}}
\rnc{\O}{\operatorname{O}}

\nc{\vtx}[2]{\fill (#1,#2)circle(.2);}
\nc{\lvtx}[2]{\fill (#1,0)circle(.2);}
\nc{\uvtx}[2]{\fill (#1,1.5)circle(.2);}

\nc{\Eq}{\mathfrak{Eq}}
\nc{\Gau}{\Ga^\wedge} \nc{\Gal}{\Ga_\vee}
\nc{\Lamu}{\Lam^\wedge} \nc{\Laml}{\Lam_\vee}
\nc{\bX}{{\bf X}}
\nc{\bY}{{\bf Y}}
\nc{\ds}{\displaystyle}

\nc{\uvert}[1]{\fill (#1,1.5)circle(.2);}
\nc{\uuvert}[1]{\fill (#1,3)circle(.2);}
\nc{\uuuvert}[1]{\fill (#1,4.5)circle(.2);}
\rnc{\lvert}[1]{\fill (#1,0)circle(.2);}
\nc{\overt}[1]{\fill (#1,0)circle(.1);}
\nc{\overtl}[3]{\node[vertex] (#3) at (#1,0) {  {\tiny $#2$} };}
\nc{\cv}[2]{\draw(#1,1.5) to [out=270,in=90] (#2,0);}
\nc{\cvs}[2]{\draw(#1,1.5) to [out=270+30,in=90+30] (#2,0);}
\nc{\ucv}[2]{\draw(#1,3) to [out=270,in=90] (#2,1.5);}
\nc{\uucv}[2]{\draw(#1,4.5) to [out=270,in=90] (#2,3);}
\nc{\textpartn}[1]{{\lower0.45 ex\hbox{\begin{tikzpicture}[xscale=.2,yscale=0.2] #1 \end{tikzpicture}}}}
\nc{\textpartnx}[2]{{\lower1.0 ex\hbox{\begin{tikzpicture}[xscale=.3,yscale=0.3]
\foreach \x in {1,...,#1}
{ \uvert{\x} \lvert{\x} }
#2 \end{tikzpicture}}}}
\nc{\disppartnx}[2]{{\lower1.0 ex\hbox{\begin{tikzpicture}[scale=0.3]
\foreach \x in {1,...,#1}
{ \uvert{\x} \lvert{\x} }
#2 \end{tikzpicture}}}}
\nc{\disppartnxd}[2]{{\lower2.1 ex\hbox{\begin{tikzpicture}[scale=0.3]
\foreach \x in {1,...,#1}
{ \uuvert{\x} \uvert{\x} \lvert{\x} }
#2 \end{tikzpicture}}}}
\nc{\disppartnxdn}[2]{{\lower2.1 ex\hbox{\begin{tikzpicture}[scale=0.3]
\foreach \x in {1,...,#1}
{ \uuvert{\x} \lvert{\x} }
#2 \end{tikzpicture}}}}
\nc{\disppartnxdd}[2]{{\lower3.6 ex\hbox{\begin{tikzpicture}[scale=0.3]
\foreach \x in {1,...,#1}
{ \uuuvert{\x} \uuvert{\x} \uvert{\x} \lvert{\x} }
#2 \end{tikzpicture}}}}

\nc{\dispgax}[2]{{\lower0.0 ex\hbox{\begin{tikzpicture}[scale=0.3]
#2
\foreach \x in {1,...,#1}
{\lvert{\x} }
 \end{tikzpicture}}}}
\nc{\textgax}[2]{{\lower0.4 ex\hbox{\begin{tikzpicture}[scale=0.3]
#2
\foreach \x in {1,...,#1}
{\lvert{\x} }
 \end{tikzpicture}}}}
\nc{\textlinegraph}[2]{{\raise#1 ex\hbox{\begin{tikzpicture}[scale=0.8]
#2
 \end{tikzpicture}}}}
\nc{\textlinegraphl}[2]{{\raise#1 ex\hbox{\begin{tikzpicture}[scale=0.8]
\tikzstyle{vertex}=[circle,draw=black, fill=white, inner sep = 0.07cm]
#2
 \end{tikzpicture}}}}
\nc{\displinegraph}[1]{{\lower0.0 ex\hbox{\begin{tikzpicture}[scale=0.6]
#1
 \end{tikzpicture}}}}

\nc{\disppartnthreeone}[1]{{\lower1.0 ex\hbox{\begin{tikzpicture}[scale=0.3]
\foreach \x in {1,2,3,5,6}
{ \uvert{\x} }
\foreach \x in {1,2,4,5,6}
{ \lvert{\x} }
\draw[dotted] (3.5,1.5)--(4.5,1.5);
\draw[dotted] (2.5,0)--(3.5,0);
#1 \end{tikzpicture}}}}

\nc{\partn}[4]{\left( \begin{array}{c|c} 
#1 \ & \ #3 \ \ \\ \cline{2-2}
#2 \ & \ #4 \ \
\end{array} \!\!\! \right)}
\nc{\partnlong}[6]{\partn{#1}{#2}{#3,\ #4}{#5,\ #6}} 
\nc{\partnsh}[2]{\left( \begin{array}{c} 
#1 \\
#2
\end{array} \right)}
\nc{\partncodefz}[3]{\partn{#1}{#2}{#3}{\emptyset}}
\nc{\partndefz}[3]{{\partn{#1}{#2}{\emptyset}{#3}}}
\nc{\partnlast}[2]{\left( \begin{array}{c|c}
#1 \ &  \ #2 \\
#1 \ &  \ #2
\end{array} \right)}

\nc{\uuarcx}[3]{\draw(#1,3)arc(180:270:#3) (#1+#3,3-#3)--(#2-#3,3-#3) (#2-#3,3-#3) arc(270:360:#3);}
\nc{\uuarc}[2]{\uuarcx{#1}{#2}{.4}}
\nc{\uuuarcx}[3]{\draw(#1,4.5)arc(180:270:#3) (#1+#3,4.5-#3)--(#2-#3,4.5-#3) (#2-#3,4.5-#3) arc(270:360:#3);}
\nc{\uuuarc}[2]{\uuuarcx{#1}{#2}{.4}}
\nc{\darcx}[3]{\draw(#1,0)arc(180:90:#3) (#1+#3,#3)--(#2-#3,#3) (#2-#3,#3) arc(90:0:#3);}
\nc{\darc}[2]{\darcx{#1}{#2}{.4}}
\nc{\udarcx}[3]{\draw(#1,1.5)arc(180:90:#3) (#1+#3,1.5+#3)--(#2-#3,1.5+#3) (#2-#3,1.5+#3) arc(90:0:#3);}
\nc{\udarc}[2]{\udarcx{#1}{#2}{.4}}
\nc{\uudarcx}[3]{\draw(#1,3)arc(180:90:#3) (#1+#3,3+#3)--(#2-#3,3+#3) (#2-#3,3+#3) arc(90:0:#3);}
\nc{\uudarc}[2]{\uudarcx{#1}{#2}{.4}}
\nc{\uarcx}[3]{\draw(#1,1.5)arc(180:270:#3) (#1+#3,1.5-#3)--(#2-#3,1.5-#3) (#2-#3,1.5-#3) arc(270:360:#3);}
\nc{\uarc}[2]{\uarcx{#1}{#2}{.4}}
\nc{\darcxhalf}[3]{\draw(#1,0)arc(180:90:#3) (#1+#3,#3)--(#2,#3) ;}
\nc{\darchalf}[2]{\darcxhalf{#1}{#2}{.4}}
\nc{\uarcxhalf}[3]{\draw(#1,1.5)arc(180:270:#3) (#1+#3,1.5-#3)--(#2,1.5-#3) ;}
\nc{\uarchalf}[2]{\uarcxhalf{#1}{#2}{.4}}
\nc{\uarcxhalfr}[3]{\draw (#1+#3,1.5-#3)--(#2-#3,1.5-#3) (#2-#3,1.5-#3) arc(270:360:#3);}
\nc{\uarchalfr}[2]{\uarcxhalfr{#1}{#2}{.4}}

\nc{\bdarcx}[3]{\draw[blue](#1,0)arc(180:90:#3) (#1+#3,#3)--(#2-#3,#3) (#2-#3,#3) arc(90:0:#3);}
\nc{\bdarc}[2]{\darcx{#1}{#2}{.4}}
\nc{\rduarcx}[3]{\draw[red](#1,0)arc(180:270:#3) (#1+#3,0-#3)--(#2-#3,0-#3) (#2-#3,0-#3) arc(270:360:#3);}
\nc{\rduarc}[2]{\uarcx{#1}{#2}{.4}}
\nc{\duarcx}[3]{\draw(#1,0)arc(180:270:#3) (#1+#3,0-#3)--(#2-#3,0-#3) (#2-#3,0-#3) arc(270:360:#3);}
\nc{\duarc}[2]{\uarcx{#1}{#2}{.4}}

\nc{\uv}[1]{\fill (#1,2)circle(.1);}
\nc{\uvw}[1]{\fill[white] (#1,2)circle(.1);}
\nc{\huv}[1]{\fill (#1,1)circle(.1);}
\nc{\lv}[1]{\fill (#1,0)circle(.1);}
\nc{\llv}[1]{\fill (#1,-2)circle(.1);}
\nc{\arcup}[2]{
\draw(#1,2)arc(180:270:.4) (#1+.4,1.6)--(#2-.4,1.6) (#2-.4,1.6) arc(270:360:.4);
}
\nc{\harcup}[2]{
\draw(#1,1)arc(180:270:.4) (#1+.4,.6)--(#2-.4,.6) (#2-.4,.6) arc(270:360:.4);
}
\nc{\arcdn}[2]{
\draw(#1,0)arc(180:90:.4) (#1+.4,.4)--(#2-.4,.4) (#2-.4,.4) arc(90:0:.4);
}
\nc{\cve}[2]{
\draw(#1,2) to [out=270,in=90] (#2,0);
}
\nc{\hcve}[2]{
\draw(#1,1) to [out=270,in=90] (#2,0);
}
\nc{\catarc}[3]{
\draw(#1,2)arc(180:270:#3) (#1+#3,2-#3)--(#2-#3,2-#3) (#2-#3,2-#3) arc(270:360:#3);
}

\nc{\arcr}[2]{
\draw[red](#1,0)arc(180:90:.4) (#1+.4,.4)--(#2-.4,.4) (#2-.4,.4) arc(90:0:.4);
}
\nc{\arcb}[2]{
\draw[blue](#1,2-2)arc(180:270:.4) (#1+.4,1.6-2)--(#2-.4,1.6-2) (#2-.4,1.6-2) arc(270:360:.4);
}
\nc{\loopr}[1]{
\draw[blue](#1,-2) edge [out=130,in=50,loop] ();
}
\nc{\loopb}[1]{
\draw[red](#1,-2) edge [out=180+130,in=180+50,loop] ();
}
\nc{\redto}[2]{\draw[red](#1,0)--(#2,0);}
\nc{\bluto}[2]{\draw[blue](#1,0)--(#2,0);}
\nc{\dotto}[2]{\draw[dotted](#1,0)--(#2,0);}
\nc{\lloopr}[2]{\draw[red](#1,0)arc(0:360:#2);}
\nc{\lloopb}[2]{\draw[blue](#1,0)arc(0:360:#2);}
\nc{\rloopr}[2]{\draw[red](#1,0)arc(-180:180:#2);}
\nc{\rloopb}[2]{\draw[blue](#1,0)arc(-180:180:#2);}
\nc{\uloopr}[2]{\draw[red](#1,0)arc(-270:270:#2);}
\nc{\uloopb}[2]{\draw[blue](#1,0)arc(-270:270:#2);}
\nc{\dloopr}[2]{\draw[red](#1,0)arc(-90:270:#2);}
\nc{\dloopb}[2]{\draw[blue](#1,0)arc(-90:270:#2);}
\nc{\llloopr}[2]{\draw[red](#1,0-2)arc(0:360:#2);}
\nc{\llloopb}[2]{\draw[blue](#1,0-2)arc(0:360:#2);}
\nc{\lrloopr}[2]{\draw[red](#1,0-2)arc(-180:180:#2);}
\nc{\lrloopb}[2]{\draw[blue](#1,0-2)arc(-180:180:#2);}
\nc{\ldloopr}[2]{\draw[red](#1,0-2)arc(-270:270:#2);}
\nc{\ldloopb}[2]{\draw[blue](#1,0-2)arc(-270:270:#2);}
\nc{\luloopr}[2]{\draw[red](#1,0-2)arc(-90:270:#2);}
\nc{\luloopb}[2]{\draw[blue](#1,0-2)arc(-90:270:#2);}

\nc{\larcb}[2]{
\draw[blue](#1,0-2)arc(180:90:.4) (#1+.4,.4-2)--(#2-.4,.4-2) (#2-.4,.4-2) arc(90:0:.4);
}
\nc{\larcr}[2]{
\draw[red](#1,2-2-2)arc(180:270:.4) (#1+.4,1.6-2-2)--(#2-.4,1.6-2-2) (#2-.4,1.6-2-2) arc(270:360:.4);
}

\rnc{\H}{\mathscr H}
\rnc{\L}{\mathscr L}
\nc{\R}{\mathscr R}
\nc{\D}{\mathscr D}
\nc{\J}{\mathscr D}

\nc{\ssim}{\mathrel{\raise0.25 ex\hbox{\oalign{$\approx$\crcr\noalign{\kern-0.84 ex}$\approx$}}}}
\nc{\POI}{\mathcal{POI}}
\nc{\wb}{\overline{w}}
\nc{\ub}{\overline{u}}
\nc{\vb}{\overline{v}}
\nc{\fb}{\overline{f}}
\rnc{\sb}{\overline{s}}
\nc{\XR}{\pres{X}{R\,}}
\nc{\YQ}{\pres{Y}{Q}}
\nc{\ZP}{\pres{Z}{P\,}}
\nc{\XRone}{\pres{X_1}{R_1}}
\nc{\XRtwo}{\pres{X_2}{R_2}}
\nc{\XRthree}{\pres{X_1\cup X_2}{R_1\cup R_2\cup R_3}}
\nc{\er}{\eqref}
\nc{\larr}{\mathrel{\hspace{-0.35 ex}>\hspace{-1.1ex}-}\hspace{-0.35 ex}}
\nc{\rarr}{\mathrel{\hspace{-0.35 ex}-\hspace{-0.5ex}-\hspace{-2.3ex}>\hspace{-0.35 ex}}}
\nc{\lrarr}{\mathrel{\hspace{-0.35 ex}>\hspace{-1.1ex}-\hspace{-0.5ex}-\hspace{-2.3ex}>\hspace{-0.35 ex}}}
\nc{\nn}{\tag*{}}
\nc{\epfal}{\tag*{$\Box$}}
\nc{\tagd}[1]{\tag*{(#1)$'$}}
\nc{\ldb}{[\![}
\nc{\rdb}{]\!]}
\nc{\sm}{\setminus}
\nc{\I}{\mathcal I}
\nc{\InSn}{\I_n\setminus\S_n}
\nc{\dom}{\operatorname{dom}^\wedge} \nc{\codom}{\operatorname{dom}_\vee}
\nc{\ojin}{1\leq j<i\leq n}
\nc{\eh}{\widehat{e}}
\nc{\wh}{\widehat{w}}
\nc{\uh}{\widehat{u}}
\nc{\vh}{\widehat{v}}
\nc{\sh}{\widehat{s}}
\nc{\fh}{\widehat{f}}
\nc{\textres}[1]{\text{\emph{#1}}}
\nc{\aand}{\emph{\ and \ }}
\nc{\iif}{\emph{\ if \ }}
\nc{\textlarr}{\ \larr\ }
\nc{\textrarr}{\ \rarr\ }
\nc{\textlrarr}{\ \lrarr\ }

\nc{\comma}{,\ }

\nc{\COMMA}{,\quad}
\nc{\TnSn}{\T_n\setminus\S_n}
\nc{\TXSX}{\T_X\setminus\S_X}
\rnc{\S}{\mathcal S}

\nc{\T}{\mathcal T}
\nc{\A}{\mathscr A}
\nc{\B}{\mathcal B}
\rnc{\P}{\mathcal P}
\nc{\K}{\mathcal K}
\nc{\PB}{\mathcal{PB}}
\nc{\rank}{\operatorname{rank}}

\nc{\mtt}{\!\!\!\mt\!\!\!}

\nc{\sub}{\subseteq}
\nc{\la}{\langle}
\nc{\ra}{\rangle}
\nc{\mt}{\mapsto}
\nc{\im}{\mathrm{im}}
\nc{\id}{\mathrm{id}}
\nc{\bn}{\mathbf{n}}
\nc{\bk}{\mathbf{k}}
\nc{\br}{\mathbf{r}}
\nc{\ve}{\varepsilon}
\nc{\al}{\alpha}
\nc{\be}{\beta}
\nc{\ga}{\gamma}
\nc{\Ga}{\Gamma}
\nc{\de}{\delta}
\nc{\ka}{\kappa}
\nc{\lam}{\lambda}
\nc{\Lam}{\Lambda}
\nc{\si}{\sigma}
\nc{\Si}{\Sigma}
\nc{\oijn}{1\leq i<j\leq n}

\nc{\comm}{\rightleftharpoons}
\nc{\AND}{\qquad\text{and}\qquad}

\nc{\bit}{\vspace{-3 truemm}\begin{itemize}}
\nc{\eit}{\end{itemize}\vspace{-3 truemm}}
\nc{\ben}{\vspace{-3 truemm}\begin{enumerate}}
\nc{\een}{\end{enumerate}\vspace{-3 truemm}}
\nc{\eitres}{\end{itemize}}

\nc{\set}[2]{\{ {#1} : {#2} \}}
\nc{\bigset}[2]{\big\{ {#1}: {#2} \big\}}
\nc{\Bigset}[2]{\Big\{ \,{#1}\, \,\Big|\, \,{#2}\, \Big\}}

\nc{\pres}[2]{\la {#1} \,|\, {#2} \ra}
\nc{\bigpres}[2]{\big\la {#1} \,\big|\, {#2} \big\ra}
\nc{\Bigpres}[2]{\Big\la \,{#1}\, \,\Big|\, \,{#2}\, \Big\ra}
\nc{\Biggpres}[2]{\Bigg\la {#1} \,\Bigg|\, {#2} \Bigg\ra}

\nc{\pf}{\noindent{\bf Proof.}  }
\nc{\epf}{\hfill$\Box$\bigskip}
\nc{\epfres}{\hfill$\Box$}
\nc{\pfnb}{\pf}
\nc{\epfnb}{\bigskip}
\nc{\pfthm}[1]{\bigskip \noindent{\bf Proof of Theorem \ref{#1}}\,\,  }
\nc{\pfprop}[1]{\bigskip \noindent{\bf Proof of Proposition \ref{#1}}\,\,  }
\nc{\epfreseq}{\tag*{$\Box$}}

\makeatletter
\newcommand\footnoteref[1]{\protected@xdef\@thefnmark{\ref{#1}}\@footnotemark}
\makeatother

\newtheorem{thm}[equation]{Theorem}
\newtheorem{lemma}[equation]{Lemma}
\newtheorem{cor}[equation]{Corollary}
\newtheorem{prop}[equation]{Proposition}

\theoremstyle{definition}

\newtheorem{rem}[equation]{Remark}
\newtheorem{defn}[equation]{Definition}
\newtheorem{eg}[equation]{Example}
\newtheorem{ass}[equation]{Assumption}

\title{Enumeration of idempotents in diagram semigroups and algebras~\vspace{-5ex}}
\author{}
\date{}

\maketitle
\begin{center}
{\large Igor Dolinka%
\footnote{Department of Mathematics and Informatics, University of Novi Sad, Trg Dositeja Obradovi\'ca 4, 21101 Novi Sad, Serbia. {\it Email:} {\tt dockie@dmi.uns.ac.rs}},
James East%
\footnote{Centre for Research in Mathematics, School of Computing, Engineering and Mathematics, University of Western Sydney, Locked Bag 1797, Penrith NSW 2751, Australia. {\it Email:} {\tt J.East@uws.edu.au}},
Athanasios Evangelou%
\footnote{School of Mathematics and Physics, University of Tasmania, Private Bag 37, Hobart 7001, Australia. {\it Email:} {\tt aoost@utas.edu.au}, {\tt D.FitzGerald@utas.edu.au}, {\tt nicholas.ham@utas.edu.au}},
Des FitzGerald${}^3$, }\\
{\large Nicholas Ham${}^3$,
James Hyde\footnote{School of Mathematics and Statistics, University of St Andrews, St Andrews, UK. {\it Email:} {\tt jth4@st-andrews.ac.uk}},
Nicholas Loughlin\footnote{School of Mathematics and Statistics, Newcastle University, Newcastle NE1 7RU, UK.  {\it Email:} {\tt n.j.loughlin@newcastle.ac.uk}}}\\
~\\
\end{center}

\begin{abstract}
We give a characterisation of the idempotents of the partition monoid, and use this to enumerate the idempotents in the finite partition, Brauer and partial Brauer monoids, giving several formulae and recursions for the number of idempotents in each monoid as well as various $\R$-, $\L$- and $\D$-classes.  We also apply our results to determine the number of idempotent basis elements in the finite dimensional partition, Brauer and partial Brauer algebras.


\end{abstract}

\section{Introduction}\label{sect:intro}

There are many compelling reasons to study diagram algebras and semigroups.  Besides their intrinsic appeal, they appear as key objects in several diverse areas of mathematics, from statistical mechanics to the representation theory of algebraic groups, often touching upon major combinatorial themes. In this introduction we seek to show the value of this study, though we can give only a superficial impression of all the connections that exist, with a particular emphasis on the types of problems we investigate, and we make no attempt to give an exhaustive description of an area that is exquisitely vast.

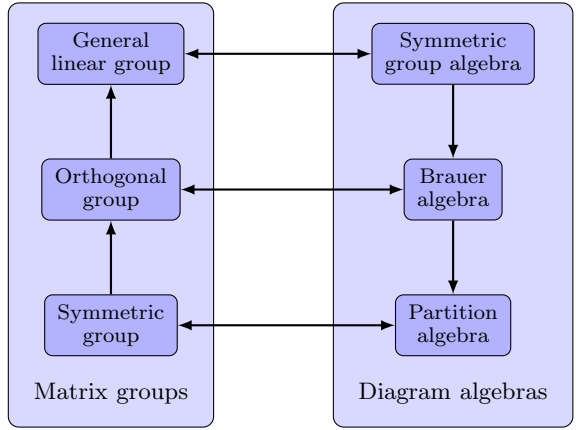
\begin{figure}[ht]
\begin{center}
\begin{tikzpicture}[scale=0.9]
\draw [fill=blue!15,rounded corners] (-1.5,-1.5)--(1.5,-1.5)--(1.5,4.75)--(-1.5,4.75)--cycle;
\draw (0,-1) node {{\footnotesize Matrix groups}};
\draw [fill=blue!15,rounded corners] (-1.75+5,-1.5)--(1.75+5,-1.5)--(1.75+5,4.75)--(-1.75+5,4.75)--cycle;
\draw (0+5,-1) node {{\footnotesize Diagram algebras}};
\path
(0,4) node(c) [rectangle,rounded corners=3pt,draw,fill=blue!30] {${\text{General} \atop \text{linear group}}$}
(0,2) node(b) [rectangle,rounded corners=3pt,draw,fill=blue!30] {${\text{Orthogonal} \atop \text{group}}$}
(0,0) node(a) [rectangle,rounded corners=3pt,draw,fill=blue!30] {${\text{Symmetric} \atop \text{group}}$}
(5,4) node(f) [rectangle,rounded corners=3pt,draw,fill=blue!30] {${\text{Symmetric} \atop \text{group algebra}}$}
(5,2) node(e) [rectangle,rounded corners=3pt,draw,fill=blue!30] {${\text{Brauer} \atop \text{algebra}}$}
(5,0) node(d) [rectangle,rounded corners=3pt,draw,fill=blue!30] {${\text{Partition} \atop \text{algebra}}$}
;
\draw[thick,{latex}-{latex}] (c) -- (f);
\draw[thick,{latex}-{latex}] (b) -- (e);
\draw[thick,{latex}-{latex}] (a) -- (d);
\draw[thick,-{latex}] (a) -- (b);
\draw[thick,-{latex}] (b) -- (c);
\draw[thick,-{latex}] (f) -- (e);
\draw[thick,-{latex}] (e) -- (d);
\end{tikzpicture}
    \caption{Schur-Weyl duality between matrix groups and diagram algebras.}
    \label{diagram_algebras}
   \end{center}
 \end{figure}

In 1927, Issai Schur \cite{Schur} provided a vital link between permutation and matrix representations.  This connection, now known as \emph{Schur-Weyl duality}, shows that the \emph{general linear group} $\GL_n(\C)$ (consisting of all invertible $n\times n$ matrices over the complex field $\C$) and the complex group algebra $\C[\S_k]$ of the \emph{symmetric group} $\S_k$ (consisting of all permutations of a $k$-element set) have commuting actions on $k$-fold tensor space $(\C^n)^{\otimes k}$, and that the irreducible components of these actions are intricately intertwined.
In 1937, Richard Brauer \cite{Brauer1937} showed that an analogous duality holds between the \emph{orthogonal group} $\operatorname{O}_n(\C)\sub\GL_n(\C)$ and the so-called \emph{Brauer algebra} ${\C^{\xi}[\B_k]\supseteq\C[\S_k]}$.
At the end of the 20th century, the \emph{partition algebras} $\C^\xi[\P_k]\supseteq\C^\xi[\B_k]$ were introduced by Paul Martin \cite{Martin1994} in the context of Potts models in statistical mechanics.  Martin later showed~\cite{Martin1998} that the partition algebras are in a kind of Schur-Weyl duality with the symmetric group $\S_n\sub\O_n(\C)$ (in its disguise as the group of all $n\times n$ permutation matrices).  Figure~\ref{diagram_algebras} shows the relationships between the various algebraic structures; vertical arrows indicate containment of algebras or groups and horizontal arrows indicate relationships between dual algebras and groups.  It should be noted that there are several other Schur-Weyl dualities; for example, between the \emph{partial Brauer algebra} $\C^\xi[\PB_k]$ and $\O_n(\C)$ \cite{MarMaz2,Haldel,GW1995,Maz1998}, and between (the semigroup algebras of) the symmetric and dual symmetric inverse semigroups \cite{KudMaz2008}.

The traditional approach to studying the above algebras, collectively referred to as \emph{diagram algebras} since they have bases indexed by certain diagrams, has been via representation theory \cite{HR,Halverson2001,Martin1996}.
But recent investigations \cite{ACHLV,FitzGerald2006,Aui2012,Aui2013,ADV,JEpme,JEgrpm,JEpnsn,JEipms,EF,EastGray,Maltcev2007,KudMaz2007,KudMaz2006,KMM2006, Maz2002,Maz1998,Larsson} have taken a more direct approach, probing the so-called \emph{partition monoids} (and other \emph{diagram semigroups}) with the tools of semigroup theory, asking and answering the same kinds of questions of partition monoids as one would of any other kind of interesting semigroup, and thereby shedding new light on the internal structure of the diagram algebras.
There are two main reasons this approach has been so successful.  The first is that the partition monoids naturally embed many important transformation semigroups on the same base set; these include the full (but not partial) transformation semigroups and the symmetric and dual symmetric inverse semigroups, allowing knowledge of these semigroups to lead to new information about the partition monoids.  (See also \cite{MarMaz} where a larger semigroup is defined that contains all of the above semigroups and more.)  The second reason is that the partition algebras have natural bases consisting of diagrams (see below for the precise definitions), with the product of two basis elements always being a scalar multiple of another.  Using this observation, Wilcox \cite{Wilcox} realised the partition algebras as \emph{twisted semigroup algebras} of the partition monoids, allowing the \emph{cellularity} of the algebras to be deduced from structural information about the associated monoid.  Cellular algebras were introduced by Graham and Lehrer \cite{GL} and provide a unified framework for studying several important classes of algebras, allowing one to obtain a great deal of information about the representation theory of the algebra; see \cite{JEcais} for the original study of cellular semigroup algebras and also \cite{Guo2009} for some recent developments.  
The elements $C_{st}^\lam$ of the cellular bases of the diagram algebras studied in \cite{Wilcox} are all sums over elements from certain ``$\H$-classes'' in a corresponding diagram semigroup.  Of importance to the cellular structure of the algebra is whether a product $C_{st}^\lam\cdot C_{uv}^\lam$ ``moves down'' in the algebra, and this is governed by the location of idempotents within the ``$\D$-class'' containing the elements involved in the sums defining $C_{st}^\lam$ and $C_{uv}^\lam$.
The twisted semigroup algebra structure has also been useful in the derivation of presentations by generators and relations \cite{JEgrpm,JEpnsn}.
But the benefits of the relationship do not only flow from semigroup theory to diagram algebras.  Indeed, the partition monoids and other kinds of diagram semigroups have played vital roles in solving outstanding problems in semigroup theory itself, especially, so far, in the context of pseudovarieties of finite semigroups \cite{ACHLV,Aui2012,Aui2013,ADV} and embeddings in regular $*$-semigroups \cite{JEpme,EF}.

It has long been recognised that the \emph{biordered set} of idempotents $E(S)=\set{x\in S}{x^2=x}$ of a semigroup $S$ often provides a great deal of useful information about the structure of the semigroup itself.  In some cases, $E(S)$ is a subsemigroup of $S$ (as in inverse semigroups, for example), but this is not generally the case.  However, the subsemigroup generated by the idempotents of a semigroup is typically a very interesting object with a rich combinatorial structure.  In many examples of finite semigroups, this subsemigroup coincides with the \emph{singular ideal}, the set of non-invertible elements \cite{Fountain1992,Howie1966,Putcha}; this is also the case with the partition and Brauer monoids \cite{EF,JEpnsn,Maltcev2007}.  Several studies have considered (minimal idempotent) generating sets of these singular ideals as well as more general ideals; see \cite{EastGray} and references therein.  Another reason idempotent generated semigroups have received considerable attention in the literature is that they possess a universal property: every semigroup embeds into an idempotent generated semigroup \cite{Howie1966} (indeed, in an idempotent generated regular $*$-semigroup~\cite{EF}). There has also been a recent resurgence of interest in the so-called \emph{free idempotent generated semigroups} (see \cite{Gray2012(2),DolinkaGray} and references therein) although, to the authors' knowledge, very little is currently known about the free idempotent generated semigroups arising from of the diagram semigroups we consider; we hope the current work will help with the pursuit of this knowledge.

Interestingly, although much is known \cite{EF,JEpnsn,EastGray,Maltcev2007} about the semigroups generated by the idempotents of certain diagram semigroups, the idempotents themselves have so far evaded classification and enumeration, apart from the case of the Brauer monoid $\B_n$ (see \cite{Larsson}, where a different approach to ours leads to sums over set partitions).  This stands in stark contrast to many other natural families of semigroup; for example, the idempotents of the symmetric inverse monoid $\I_X$ are the restrictions of the identity map, while the idempotents of the full transformation semigroup~$\T_X$ are the transformations that map their image identically.  These descriptions allow for easy enumeration; for example, $|E(\I_n)|=2^n$, and $|E(\T_n)|=\sum_{k=1}^n{n\choose k}k^{n-k}$.  It is the goal of this article to rectify the situation for several classes of diagram semigroups; specifically, the partition, Brauer and partial Brauer monoids $\P_n,\B_n,\PB_n$ (though much of what we say will also apply to various transformation semigroups such as $\I_n$ and~$\T_n$).  For each semigroup, we completely describe the idempotents, and we give several formulae and recursions for the number of idempotents in the semigroup as well as in various Green's classes and ideals.  We also give formulae for the number of idempotent basis elements in the corresponding diagram algebras; these depend on whether the constant that determines the twisting is a root of unity.  Our approach is combinatorial in nature, and our results depend on certain equivalence relations and graphs associated to a partition.  
Because Sloane's Online Encyclopedia of Integer Sequences (OEIS) \cite{OEIS} is an important resource in many areas of discrete mathematics, we record the sequences that result from our study.
We remark that our approach does not work for the so-called \emph{Jones monoid} $\mathcal J_n\sub\B_n$ (also sometimes called the \emph{Temperley-Lieb monoid} and denoted $\mathcal{T}\!\mathcal{L}_n$), which consists of all \emph{planar} Brauer diagrams; values of $|E(\mathcal J_n)|$ up to $n=19$ have been calculated by James Mitchell, using the {\sf Semigroups} package in {\sf GAP} \cite{GAP}, and may be found in Sequence A225798 on the OEIS \cite{OEIS}.

The article is organised as follows.  In Section \ref{sect:prelim}, we define the diagram semigroups we will be studying, and we state and prove some of the basic properties we will need.  The characterisation of the idempotents is given in Section~\ref{sect:char}, with the main result being Theorem~\ref{char}.  In Section~\ref{sect:enum}, we enumerate the idempotents, first giving general results (Theorems~\ref{eKn},~\ref{eDrKn},~\ref{eDrKn_rec}) and then applying these to the partition, Brauer and partial Brauer monoids in Sections~\ref{subsect:Pn},~\ref{subsect:Bn} and~\ref{subsect:PBn}.  We describe an alternative approach to the enumeration of the idempotents in the Brauer and partial Brauer monoids in Section~\ref{sect:diff} (see Theorems~\ref{anr} and~\ref{anrt}).  In Section~\ref{sect:algebras}, we classify and enumerate the idempotent basis elements in the partition, Brauer and partial Brauer algebras (see especially Theorems~\ref{char2},~\ref{eKnxi=0},~\ref{eKnxi>0}).  Finally, in Section~\ref{sect:values}, we give several tables of calculated values.
The reader is referred to the monographs \cite{Hig,Howie} for background on semigroups.  


\section{Preliminaries}\label{sect:prelim}

Let $X$ be a set, and $X'$ a disjoint set in one-one correspondence with $X$ via a mapping $X\to X':x\mt x'$.  If $A\sub X$ we will write $A'=\set{a'}{a\in A}$. A \emph{partition on $X$} is a collection of pairwise disjoint non-empty subsets of $X\cup X'$ whose union is $X\cup X'$; these subsets are called the \emph{blocks} of the partition.  The \emph{partition monoid} on $X$ is the set $\P_X$ of all such partitions, with a natural binary operation defined below. 
When $n\in\mathbb N=\{0,1,2,\ldots\}$ is a natural number and $X=\bn=\{1,\ldots,n\}$, we will write $\P_X=\P_n$.  Note that $\P_0=\P_\emptyset = \{\emptyset\}$ has a single element; namely, the empty partition, which we denote by $\emptyset$.

A partition may be represented as a graph on the vertex set $X\cup X'$; edges are included so that the connected components of the graph correspond to the blocks of the partition.  Of course such a graphical representation is not unique, but we regard two such graphs as equivalent if they have the same connected components, and we typically identify a partition with any graph representing it.  We think of the vertices from~$X$ (resp.~$X'$) as being the \emph{upper vertices} (resp.~\emph{lower vertices}).  For example, the partition
$\al=\big\{ \{1,4\},\{2,3,4',5'\},\{5,6\},\{1',3',6'\},\{2'\}\big\}\in\P_6$
is represented by the graph
$\al=\textpartnx6{
\uarcx14{.4}
\uarcx23{.2}
\uarcx56{.2}
\darc13
\darc36
\darcx45{.2}
\cv34
}$.

In order to describe the product alluded to above, let $\al,\be\in\P_X$. Consider now a third set~$X''$, disjoint from both $X$ and $X'$, and in bijection with $X$ via $x\mt x''$.  Let $\alb$ be the graph obtained from (a graph representing) $\al$ simply by changing the label of each lower vertex $x'$ to~$x''$.  Similarly, let $\beb$ be the graph obtained from $\be$ by changing the label of each upper vertex~$x$ to~$x''$.  Consider now the graph $\Ga(\al,\be)$ on the vertex set~$X\cup X'\cup X''$ obtained by joining $\alb$ and~$\beb$ together so that each lower vertex $x''$ of $\alb$ is identified with the corresponding upper vertex $x''$ of $\beb$.  Note that $\Ga(\al,\be)$, which we call the \emph{product graph} of $\al$ and $\be$, may contain multiple edges.  We define $\al\be\in\P_X$ to be the partition that satisfies the property that $x,y\in X\cup X'$ belong to the same block of $\al\be$ if and only if there is a path from $x$ to $y$ in $\Ga(\al,\be)$.
An example calculation (with $X$ finite) is given in Figure~\ref{fig:multinP6}.
\begin{figure}[ht]
   \begin{center}
$$
\disppartnxdd6{
\uuuarcx14{.4}
\uuuarcx23{.2}
\uuuarcx56{.2}
\uudarc13
\uudarc36
\uudarcx45{.2}
\uucv34
\uarcx13{.4}
\uarcx24{.6}
\darcx45{.2}
\darcx56{.2}
\cv43
\cv55
\draw(0.7,3.75)node[left]{{\small $\al=$}};
\draw(0.7,.75)node[left]{{\small $\be=$}};
}
\quad\longrightarrow\quad
\disppartnxd6{
\uuarcx14{.4}
\uuarcx23{.2}
\uuarcx56{.2}
\udarc13
\udarc36
\udarcx45{.2}
\ucv34
\uarcx13{.4}
\uarcx24{.6}
\darcx45{.2}
\darcx56{.2}
\cv43
\cv55
}
\quad\longrightarrow\quad
\disppartnx6{
\uarcx14{.4}
\uarcx23{.2}
\uarcx56{.2}
\darcx34{.2}
\darcx45{.2}
\darcx56{.2}
\cv33
\draw(6.3,0.75)node[right]{{\small $=\al\be$}};
}
$$
    \caption{Two partitions $\al,\be\in\P_6$ (left), their product $\al\be\in\P_6$ (right), and the product graph $\Ga(\al,\be)$ (centre).}
    \label{fig:multinP6}
   \end{center}
 \end{figure}

We now define subsets
\begin{align*}
\PB_X &= \set{\al\in\P_X}{\text{each block of $\al$ has size at most $2$}} \\
\B_X &= \set{\al\in\P_X}{\text{each block of $\al$ has size $2$}}.
\end{align*}
We note that $\PB_X$ is a submonoid of $\P_X$ for any set $X$, while $\B_X$ is a submonoid if and only if~$X$ is finite.  For example, taking $X=\mathbb N=\{0,1,2,3,\ldots\}$, the partitions $\al,\be$ pictured in Figure~\ref{fig:BN} both belong to~$\B_\N$, while the product $\al\be$ (also pictured in Figure \ref{fig:BN}) does not.  We call $\PB_X$ the \emph{partial Brauer monoid} and (in the case that $X$ is finite) $\B_X$ the \emph{Brauer monoid} on $X$.  Again, if $n\in\mathbb N$ and $X=\bn$, we write $\PB_n$ and $\B_n$ for $\PB_X$ and $\B_X$, noting that $\B_0=\PB_0=\P_0=\{\emptyset\}$.

\begin{figure}[ht]
   \begin{center}
$$
\disppartnxdd7{
\draw[dotted](8,4.5)--(10,4.5);
\draw[dotted](8,3)--(10,3);
\draw[dotted](8,1.5)--(10,1.5);
\draw[dotted](8,0)--(10,0);
\uuuarc23
\uuuarc45
\uuuarc67
\uudarc23
\uudarc45
\uudarc67
\uucv11
\uarc12
\uarc34
\uarc56
\darc12
\darc34
\darc56
\darchalf7{7.5}
\uarchalf7{7.5}
\draw(0.7,3.75)node[left]{{\small $\al=$}};
\draw(0.7,.75)node[left]{{\small $\be=$}};
}
\quad\longrightarrow\quad
\disppartnxd7{
\draw[dotted](8,3)--(10,3);
\draw[dotted](8,1.5)--(10,1.5);
\draw[dotted](8,0)--(10,0);
\uuarc23
\uuarc45
\uuarc67
\udarc23
\udarc45
\udarc67
\ucv11
\uarc12
\uarc34
\uarc56
\darc12
\darc34
\darc56
\darchalf7{7.5}
\uarchalf7{7.5}
}
\quad\longrightarrow\quad
\disppartnx7{
\draw[dotted](8,1.5)--(10,1.5);
\draw[dotted](8,0)--(10,0);
\uarc23
\uarc45
\uarc67
\darc12
\darc34
\darc56
\darchalf7{7.5}
\draw(10.3,0.75)node[right]{{\small $=\al\be$}};
}
$$
    \caption{Two partitions $\al,\be\in\B_\N$ (left), their product $\al\be\not\in\B_\N$ (right), and the product graph $\Ga(\al,\be)$ (centre).}
    \label{fig:BN}
   \end{center}
 \end{figure}

We now introduce some notation and terminology that we will use throughout our study.  Let $\al\in\P_X$. A block $A$ of $\al$ is called a \emph{transversal block} if $A\cap X\not=\emptyset\not= A\cap X'$, or otherwise an \emph{upper} (resp.~\emph{lower}) \emph{non-transversal block} if $A\cap X'=\emptyset$ (resp.~$A\cap X=\emptyset$).  The \emph{rank} of $\al$, denoted $\rank(\al)$, is equal to the number of transversal blocks of $\al$.  For $x\in X\cup X'$, let $[x]_\al$ denote the block of $\al$ containing $x$.  We define the \emph{upper} and \emph{lower domains} of $\al$ to be the sets
\begin{align*}
\dom(\al) = \bigset{ x\in X } { [x]_\al\cap X'\not=\emptyset} \qquad&\text{and}\qquad \codom(\al) = \bigset{ x\in X } { [x']_\al\cap X\not=\emptyset}. \\
\intertext{We also define the \emph{upper} and \emph{lower kernels} of $\al$ to be the equivalences}
\keru(\al) = \bigset{(x,y)\in X\times X}{[x]_\al=[y]_\al} \qquad&\text{and}\qquad
\kerl(\al) = \bigset{(x,y)\in X\times X}{[x']_\al=[y']_\al}.
\end{align*}
(The upper and lower domains and the upper and lower kernels have been called the \emph{domain}, \emph{codomain}, \emph{kernel} and \emph{cokernel} (respectively) in other works \cite{EastGray,JEgrpm,JEpnsn,EF}, but there should be no confusion.)
To illustrate these definitions, consider the partition $\al=\textpartnx6{
\uarcx14{.4}
\uarcx23{.2}
\uarcx56{.2}
\darc13
\darc36
\darcx45{.2}
\cv34
}$
from $\P_6$.  Then $\rank(\al)=1$, $\dom(\al)=\{2,3\}$, $\codom(\al)=\{4,5\}$, and $\al$ has upper kernel-classes~$\{1,4\}$, $\{2,3\}$, $\{5,6\}$, and lower kernel-classes $\{1,3,6\}$, $\{2\}$, $\{4,5\}$.

Recall that \emph{Green's relations} are defined on a semigroup $S$ by
\begin{gather*}
x\R y \iff x\Sone=y\Sone \COMMA x\L y \iff \Sone x=\Sone y \COMMA x\mathscr J y \iff \Sone x\Sone =\Sone y\Sone  , \\
\H=\R\cap\L \COMMA \D=\R\circ\L=\L\circ\R.
\end{gather*}
Here, $\Sone$ denotes the monoid obtained from $S$ by adjoining an identity element $1$, if necessary.
If~$S$ is finite, then $\mathscr J=\D$.  For more on Green's relations, the reader is referred to \cite{Hig,Howie}.
The next result was first proved for finite $X$ in \cite{Wilcox,Larsson}, and then in full generality in \cite{FitLau}, though the language used in those papers was different from ours; see also \cite{Maz1998} on finite (partial and full) Brauer monoids.  

\bs
\begin{thm}[{\cite[Theorem~17]{Wilcox}}]
\label{green}
For each $\alpha, \beta \in \P_X$, we have
\ben
\item[\emph{(i)}] $\alpha \R \beta$ if and only if $\dom(\alpha)=\dom(\beta)$ and $\keru(\alpha)=\keru(\beta)$\emph{;}
\item[\emph{(ii)}] $\alpha \L \beta$ if and only if $\codom(\alpha)=\codom(\beta)$ and $\kerl(\alpha)=\kerl(\beta)$\emph{;}
\item[\emph{(iii)}] $\alpha \J \beta$ if and only if $\rank(\alpha) = \rank(\beta)$. \epfres
\een
\end{thm}

Finally, we define the \emph{kernel} of $\al$ to be the join
\begin{align*}
\KER(\al) &= \keru(\al)\vee\kerl(\al).
\end{align*}
(The \emph{join} $\ve\vee\eta$ of two equivalence relations $\ve,\eta$ on $X$ is the smallest equivalence relation containing the union $\ve\cup\eta$; that is, $\ve\vee\eta$ is the transitive closure of $\ve\cup\eta$.)  The equivalence classes of $X$ with respect to $\ker(\al)$ are called the \emph{kernel-classes} of $\al$.  We call a partition $\al\in\P_X$ \emph{irreducible} if it has only one kernel-class; that is, $\al$ is irreducible if and only if $\ker(\al)=X\times X$.  Some (but not all) partitions from $\P_X$ may be built up from irreducible partitions in a way we make precise below.

The equivalences $\keru(\al),\kerl(\al),\ker(\al)$ may be visualised graphically as follows.  We define a graph $\Gau(\al)$ with vertex set $X$, and red edges drawn so that the connected components are precisely the $\keru(\al)$-classes of $X$, and we define $\Gal(\al)$ analogously but with blue edges.  (Again, there are several possible choices for $\Gau(\al)$ and $\Gal(\al)$, but we regard them all as equivalent.)  
We also define $\Ga(\al)$ to be the graph on vertex set $X$ with all the edges from both $\Gau(\al)$ and $\Gal(\al)$.  Then the connected components of $\Ga(\al)$ are precisely the kernel-classes of $\al$.  
To illustrate these ideas, consider the partitions $\al=\textpartnx6{
\uarcx14{.4}
\uarcx23{.2}
\uarcx56{.2}
\darc13
\darc36
\darcx45{.2}
\cv34
}$ and $\be=\textpartnx6{
\uarcx13{.4}
\uarcx24{.6}
\darcx45{.2}
\darcx56{.2}
\cv43
\cv55
}$
from $\P_6$.  Then $\Ga(\al)=\textgax6{
\rduarcx14{.6}
\rduarcx23{.2}
\rduarcx56{.2}
\bdarcx13{.4}
\bdarcx36{.4}
\bdarcx45{.2}
}$ and $\Ga(\be)=
\textgax6{
\rduarcx13{.4}
\rduarcx24{.6}
\bdarcx45{.2}
\bdarcx56{.2}
}$.  So $\al$ is irreducible but $\be$ is not.

It will be convenient to conclude this section with two technical results that will help simplify subsequent proofs.

\ms\smallskip
\begin{lemma}\label{tech1}
Let $\al,\be\in\P_X$ and suppose $x,y\in X$.  Then $(x,y)\in\kerl(\al)\vee\keru(\be)$ if and only if $x''$ and $y''$ are joined by a path in the product graph $\Ga(\al,\be)$.
\end{lemma}

\pf If $(x,y)\in\kerl(\al)\vee\keru(\be)$, then there is a sequence $x=x_0,x_1,\ldots,x_k=y$ such that $(x_0,x_1)\in\kerl(\al)$, $(x_1,x_2)\in\keru(\be)$, $(x_2,x_3)\in\kerl(\al)$, and so on.  Such a sequence gives rise to a path $x''=x_0''\to x_1''\to\cdots\to x_k''=y''$ in the product graph $\Ga(\al,\be)$.

Conversely, suppose $x''$ and $y''$ are joined by a path in the product graph $\Ga(\al,\be)$.  We prove that $(x,y)\in\kerl(\al)\vee\keru(\be)$ by induction on the length of a path $x''=z_0\to z_1\to\cdots\to z_t=y''$ in $\Ga(\al,\be)$.  If $t=0$, then $x=y$, and we are done, so suppose $t\geq1$.  If $z_r=w''$ for some $0<r<t$, then an induction hypothesis applied to the shorter paths $x''\to\cdots\to w''$ and $w''\to\cdots\to y''$ tells us that $(x,w)$ and $(w,y)$, and hence also $(x,y)$, belong to $\kerl(\al)\vee\keru(\be)$.  If none of $z_1,\ldots,z_{t-1}$ belong to $X''$, then they either all belong to $X$ or all to $X'$.  In the former case, it follows that $z_1,\ldots,z_{t-1},y'\in[x']_\al$, so that $(x,y)\in\kerl(\al)\sub\kerl(\al)\vee\keru(\be)$.  The other case is similar. \epf

\begin{lemma}\label{tech2}
Let $\al,\be\in\P_X$ and suppose $A\cup B'$ is a transversal block of $\al\be$.  Then for any $a\in A$ and $b\in B$, and any $c,d\in X$ with $c'\in[a]_\al$, $d\in[b']_\be$, we have $(c,d)\in\kerl(\al)\vee\keru(\be)$.
 \end{lemma}

\pf Consider a path $a=z_0\to z_1\to\cdots\to z_k=b'$ in the product graph $\Ga(\al,\be)$.
Let $0\leq r\leq k$ be the least index for which $z_r$ does not belong to $X$, and let $0\leq s\leq k$ be the greatest index for which $z_s$ does not belong to $X'$.  Then $z_r=x''$ and $z_s=y''$ for some $x,y\in X$ with $x'\in[a]_\al$ and $y\in[b']_\be$.  Since there is a path from $x''$ to $y''$ in $\Ga(\al,\be)$, Lemma \ref{tech1} tells us that $(x,y)\in\kerl(\al)\vee\keru(\be)$.  But we also have $(c,x)\in\kerl(\al)$ and $(y,d)\in\keru(\be)$.  Putting this all together gives $(c,d)\in\kerl(\al)\vee\keru(\be)$, as required. \epf

\bsa
\section{Characterisation of idempotents}\label{sect:char}

We now aim to give a characterisation of the idempotent partitions, and our first step in this direction is to describe the irreducible idempotents.  (Recall that $\al\in\P_X$ is irreducible if $\ker(\al)=X\times X$.)

\ms\smallskip
\begin{lemma}\label{tech3}
Suppose $\al\in\P_X$ is irreducible.  Then $\al$ is an idempotent if and only if $\rank(\al)\leq1$.
\end{lemma}

\pf It is clear that any partition of rank $0$ is idempotent.  Next, suppose $\rank(\al)=1$, and let the unique transversal block of $\al$ be $A\cup B'$.  Every non-transversal block of $\al$ is a block of $\al^2$, so it suffices to show that $A\cup B'$ is a block of $\al^2$.  So suppose $a\in A$ and $b\in B$.  Then there is a path from $a$ to $b'$ in (a graph representing) $\al$, so it follows that there is a path from $a$ to $b''$ and a path from $a''$ to $b'$ in the product graph $\Ga(\al,\al)$.  Since $\al$ is irreducible, $\kerl(\al)\vee\keru(\al)=\KER(\al)=X\times X$, so Lemma \ref{tech1} says that there is also a path from $b''$ to $a''$.  Putting these together, we see that there is a path from $a$ to $b'$.  This completes the proof that $\al$ is idempotent.

Now suppose $\rank(\al)\geq2$ and let $A\cup B'$ and $C\cup D'$ be distinct transversal blocks of $\al$.  Let $a,b,c,d$ be arbitrary elements of $A,B,C,D$, respectively.  Since $\al$ is irreducible, there is a path from $b''$ to $d''$ in the product graph $\Ga(\al,\al)$.  But since $a\in A$ and $b\in B$, there is a path from $a$ to $b''$ in $\Ga(\al,\al)$, and similarly there is a path from $d''$ to $c$.  Putting these together, we see that there is a path from $a$ to $c$, so that $A\cup C$ is contained in a block of $\al^2$.  But $A$ and $C$ are contained in different blocks of $\al$, so it follows that $\al$ could not be an idempotent.
\epf

We now show how idempotent partitions are built up out of irreducible ones.  Suppose $X_i$ ($i\in I$) is a family of pairwise disjoint sets, and write $X=\bigcup_{i\in I}X_i$.  We define
\[
\bigoplus_{i\in I}\P_{X_i} = \set{\al\in\P_X}{\text{each block of $\al$ is contained in $X_i\cup X_i'$ for some $i\in I$}},
\]
which is easily seen to be a submonoid of $\P_X$, and isomorphic to the direct product $\prod_{i\in I}\P_{X_i}$.  Suppose $\al\in\bigoplus_{i\in I}\P_{X_i}$.  For each $i\in I$, let $\al_i= \set{A\in\al}{A\sub X_i\cup X_i'}\in\P_{X_i} $.  We call $\al_i$ the \emph{restriction of $\al$ to $X_i$}, and we write $\al_i=\al|_{X_i}$ and $\al=\bigoplus_{i\in I}\al_i$.
We are now ready to prove the main result of this section, which gives a characterisation of the idempotent partitions.
A precursor also appears in \cite{Maz1998,Larsson} for finite (partial and full) Brauer monoids.  

\bs
\begin{thm}\label{char}
Let $\al\in\P_X$, and suppose the kernel-classes of $\al$ are $X_i\ (i\in I)$.  Then $\al$ is an idempotent if and only if the following two conditions are satisfied:
\ben
\item[\emph{(i)}] $\al\in\bigoplus_{i\in I}\P_{X_i}$, and
\item[\emph{(ii)}] the restrictions $\al|_{X_i}$ all have rank at most $1$.
\een
\end{thm}

\pf Suppose first that $\al$ is an idempotent, but that condition (i) fails.  Then there is a block $A\cup B'$ of $\al$ such that $A\sub X_i$ and $B\sub X_j$ for distinct $i,j\in I$.  Let $a\in A$ and $b\in B$.  Since $\al$ is an idempotent, $A\cup B'$ is a block of $\al^2$, and we also have $b'\in[a]_\al$ and $a\in[b']_\al$.  So Lemma~\ref{tech2} tells us that $(a,b)\in\kerl(\al)\vee\keru(\al)=\KER(\al)$.  But this contradicts the fact that $a\in X_i$ and $b\in X_j$, with $X_i$ and $X_j$ distinct kernel-classes.  Thus, (i) must hold.  It follows that $\al=\bigoplus_{i\in I}\al_i$ where $\al_i=\al|_{X_i}$ for each $i$.  Then $\bigoplus_{i\in I}\al_i=\al=\al^2=\bigoplus_{i\in I}\al_i^2$, so that each $\al_i$ is an irreducible idempotent, and (ii) now follows from Lemma \ref{tech3}.

Conversely, suppose (i) and (ii) both hold, and write $\al=\bigoplus_{i\in I}\al_i$ where $\al_i\in\P_{X_i}$ for each~$i$.  Since $\rank(\al_i)\leq1$, Lemma \ref{tech3} says that each $\al_i$ is an idempotent.  It follows that $\al^2=\bigoplus_{i\in I}\al_i^2=\bigoplus_{i\in I}\al_i=\al$.
\epf

\bsa
\section{Enumeration of idempotents}\label{sect:enum}

For a subset $\Si$ of the partition monoid $\P_X$, we write $E(\Si)=\set{\al\in \Si}{\al^2=\al}$ for the set of all idempotents from $\Si$, and we write $e(\Si)=|E(\Si)|$.  In this section, we aim to derive formulae for $e(\K_X)$ where $\K_X$ is one of $\P_X,\B_X,\PB_X$.  The infinite case is essentially trivial, but we include it for completeness.

\ms\smallskip
\begin{prop}
If $X$ is infinite, then $e(\B_X)=e(\PB_X)=e(\P_X)=2^{|X|}$.
\end{prop}

\pf Since $\B_X\sub\PB_X\sub\P_X$ and $|\P_X|=2^{|X|}$, it suffices to show that $e(\B_X)=2^{|X|}$.  Let $\A=\set{A\sub X}{|X\sm A|\geq\aleph_0}$, and let $A\in\A$.  Let $\be_A$ be any element of $\B_X$ with $\dom(\be_A)=\codom(\be_A)=A$ and such that $\{a,a'\}$ is a block of $\be_A$ for all $a\in A$.  Then $\be_A$ is clearly an idempotent.  The map $\A\to E(\B_X):A\mt\be_A$ is clearly injective, so the result follows since $|\A|=2^{|X|}$. \epf

The rest of the paper concerns the finite case so, unless stated otherwise, $X$ will denote a finite set from here on.

For a subset $\Si$ of $\P_X$, we write $C(\Si)$ for the set of all irreducible idempotents of $\Si$.  So, by Theorem~\ref{char}, $C(\Si)$ consists of all partitions $\al\in\Si$ such that $\ker(\al)=X\times X$ and $\rank(\al)\leq1$.  We will also write $c(\Si)=|C(\Si)|$.  Our next goal is to show that we may deduce the value of $e(\K_n)$ from the values of $c(\K_n)$ when $\K_n$ is one of $\P_n,\B_n,\PB_n$.

Recall that an \emph{integer partition of $n$} is a $k$-tuple $\mu=(m_1,\ldots,m_k)$ of integers, for some $k$, satisfying $m_1\geq\cdots\geq m_k\geq1$ and $m_1+\cdots+m_k=1$.  We write $\mu\vdash n$ to indicate that $\mu$ is an integer partition of $n$.  With $\mu\vdash n$ as above, we will also write $\mu=(1^{\mu_1},\ldots,n^{\mu_n})$ to indicate that, for each~$i$, exactly~$\mu_i$ of the $m_j$ are equal to $i$.
By convention, we consider $\mu=\emptyset$ to be the unique integer partition of~$0$.

Recall that a \emph{set partition of $X$} is a collection $\bX=\set{X_i}{i\in I}$ of pairwise disjoint non-empty subsets of $X$ whose union is $X$.   We will write $\bX\models X$ to indicate that $\bX$ is a set partition of $X$.  Suppose $\bX=\{X_1,\ldots,X_k\}\models\bn$.  For $i\in\bn$, write $\mu_i(\bX)$ for the cardinality of the set $\set{j\in\bk}{|X_j|=i}$, and put $\mu(\bX) = (1^{\mu_1(\bX)},\ldots,n^{\mu_n(\bX)})$, so $\mu(\bX)\vdash n$.
For $\mu=(1^{\mu_1},\ldots,n^{\mu_n})\vdash n$, we write $\pi(\mu)$ for the number of set partitions $\bX\models\bn$ such that $\mu(\bX)=\mu$.  It is  easily seen (and well-known) that
\[
\pi(\mu)
= \frac{n!}{\prod_{i=1}^n\mu_i!(i!)^{\mu_i}} .
\]
If $\al\in\P_n$ has kernel classes $X_1,\ldots,X_k$, we write $\mu(\al)=\mu(\bX)$ where $\bX=\{X_1,\ldots,X_k\}$.  Note that if $|X_1|\geq\cdots\geq|X_k|$, then $\mu(\al)=(|X_1|,\ldots,|X_k|)$ in the alternative notation for integer partitions.

\bs
\begin{thm}\label{eKn}
If $\K_n$ is one of $\P_n,\B_n,\PB_n$, then
\[
e(\K_n)
= n!\cdot \sum_{\mu\vdash n} \prod_{i=1}^n \frac{c(\K_i)^{\mu_i}}{\mu_i! (i!)^{\mu_i}} .
\]
The numbers $e(\K_n)$ satisfy the recurrence:
\[
e(\K_0)=1, \quad e(\K_n) = \sum_{m=1}^n {n-1\choose m-1} c(\K_m)e(\K_{n-m}) \quad\text{for $n\geq1$.}
\]
The values of $c(\K_n)$ are given in Propositions \ref{cPn}, \ref{cBn} and \ref{cPBn}.
\end{thm}

\pf
Fix an integer partition $\mu=(m_1,\ldots,m_k)=(1^{\mu_1},\ldots,n^{\mu_n})\vdash n$.  We count the number of idempotents $\al$ from $\K_n$ with $\mu(\al)=\mu$.  We first choose the kernel-classes $X_1,\ldots,X_k$ of $\al$, with $|X_i|=m_i$ for each $i$, which we may do in $\pi(\mu)$ ways.  For each $i$, the restriction of $\al$ to $X_i$ is an irreducible idempotent of $\K_{X_i}$, and there are precisely $c(\K_{X_i}) = c(\K_{m_i})$ of these.  So there are $c(\K_{m_1})\cdots c(\K_{m_k})=c(\K_1)^{\mu_1}\cdots c(\K_n)^{\mu_n}$ idempotents with kernel classes $X_1,\ldots,X_k$.  Multiplying by $\pi(\mu)$ and summing over all $\mu$ gives the first equality.

For the recurrence, note first that $E(\K_0)=\K_0=\{\emptyset\}$, where $\emptyset$ denotes the empty partition.  Now suppose $n\geq1$ and let $m\in\bn$.  We will count the number of idempotents $\al$ from $\K_n$ such that the kernel-class $A$ of $\al$ containing $1$ has size $m$.  We first choose the remaining $m-1$ elements of~$A$, which we may do in ${n-1\choose m-1}$ ways.  The restriction $\al|_A$ is an irreducible idempotent from $\K_A$, and may be chosen in $c(\K_A)=c(\K_m)$ ways, while the restriction $\al|_{\bn\sm A}$ is an idempotent from~$\K_{\bn\sm A}$, and may be chosen in $e(\K_{\bn\sm A})=e(\K_{n-m})$ ways. Summing over all $m\in\bn$ gives the recurrence, and completes the proof. \epf

It will also be convenient to record a result concerning the number of idempotents of a fixed rank.  If $\K_n$ is one of $\P_n,\B_n,\PB_n$ and $0\leq r\leq n$, we write
\[
\DrKn = \set{\al\in\K_n}{\rank(\al)=r}.
\]
So, by Theorem \ref{green} and the fact that $\B_n$ and $\PB_n$ are regular subsemigroups of $\P_n$, the sets $\DrKn$ are precisely the $\D$-classes of $\K_n$.  Note that $\DrBn$ is non-empty if and only if $n\equiv r\pmod 2$, as each non-transversal block of an element of $\B_n$ has size $2$.  For a subset $\Si\sub\P_n$ and an integer partition $\mu\vdash n$, we write
\[
E_\mu(\Si) = \set{\al\in E(\Si)}{\mu(\al)=\mu} \AND e_\mu(\Si) = |E_\mu(\Si)|.
\]
For a subset $\Si$ of $\P_n$, and for $r\in\{0,1\}$, let
\[
C_r(\Si) = \set{\al\in C(\Si)}{\rank(\al)=r} \AND c_r(\Si)=|C_r(\Si)|.
\]
So by Lemma \ref{tech3}, $c(\Si)=c_0(\Si)+c_1(\Si)$.

\bs
\begin{thm}\label{eDrKn}
Suppose $\K_n$ is one of $\P_n,\B_n,\PB_n$ and $0\leq r\leq n$.  Then $$e(\DrKn) = \sum_{\mu\vdash n} e_\mu(\DrKn).$$  If $\mu=(m_1,\ldots,m_k)=(1^{\mu_1},\ldots,n^{\mu_n})\vdash n$, then
\[
e_\mu(\DrKn) = 
 \frac{n!}{\prod_{i=1}^n\mu_i!(i!)^{\mu_i}}
 \sum_{A\sub\bk \atop |A|=r} \left(\prod_{i\in A}c_1(\K_{m_i})\cdot\prod_{j\in\bk\sm A}c_0(\K_{m_j})\right).
\]
The values of $c_r(\K_n)$ are given in Propositions \ref{cPn}, \ref{cBn} and \ref{cPBn}.
\end{thm}

\pf Fix $\mu=(1^{\mu_1},\ldots,n^{\mu_n})=(m_1,\ldots,m_k)\vdash n$, and suppose $\al\in E(\DrKn)$ is such that $\mu(\al)=\mu$.  We choose the kernel-classes $X_1,\ldots,X_k$ (where $|X_i|=m_i$) of $\al$ in $\pi(\mu)$ ways.  Now, $\al=\al_1\oplus\cdots\oplus\al_k$, with $\al_i\in C(\K_{X_i})$ for each $i$, and $r=\rank(\al)=\rank(\al_1)+\cdots+\rank(\al_k)$.  So we require precisely $r$ of the $\al_i$ to have rank $1$.  For each subset $A\sub\bk$ with $|A|=r$, there are $\prod_{i\in A}c_1(\K_{m_i})\cdot\prod_{j\in\bk\sm A}c_0(\K_{m_j})$ ways to choose the $\al_i$ so that $\rank(\al_i)=1$ if and only if $i\in A$.  Summing over all such $A$ gives the result.~\epf

\begin{rem}
If $\K_n$ is one of $\P_n,\B_n,\PB_n$, then the ideals of $\K_n$ are precisely the sets $$I_r(\K_n)=\bigcup_{s\leq r}D_s(\K_n).$$  (See \cite{FitzGerald2006,EastGray}.)  It follows immediately that the number of idempotents in such an ideal is given by $e(I_r(\K_n)) = \sum_{s\leq r}e(D_s(\K_n))$, so these values may be deduced from the values of $e(D_s(\K_n))$ given above.
\end{rem}

We also give a recurrence for the numbers $e(\DrKn)$.

\ms\smallskip
\begin{thm}\label{eDrKn_rec}
The numbers $e(\DrKn)$ satisfy the recurrence
\begin{align*}
e(D_n(\K_n)) &= 1 \\ 
e(D_0(\K_n)) &= \rho(\K_n)^2 \\ 
e(\DrKn) &= \sum_{m=1}^n {n-1 \choose m-1} \Big(c_0(\K_m)e(D_r(\K_{n-m})) + c_1(\K_m)e(D_{r-1}(\K_{n-m}))\Big) &&\text{if $1\leq r\leq n-1$,}
\end{align*}
where $\rho(\K_n)$ is the number of $\R$-classes in $D_0(\K_n)$; these values are given in Lemma \ref{lamKn}.  The values of $c_r(\K_n)$ are given in Propositions \ref{cPn}, \ref{cBn} and \ref{cPBn}.
\end{thm}

\pf Note that $D_n(\K_n)$ is the group of units of $\K_n$ (which is the symmetric group $\S_n$), so $e(D_n(\K_n))=1$ for all $n$.  Also, since every element of $D_0(\K_n)$ is an idempotent, it follows that $e(D_0(\K_n))=|D_0(\K_n)|=\rho(\K_n)^2$.  Now consider an element $\al\in\DrKn$ where $1\leq r\leq n-1$, and suppose the kernel-class $A$ of $\al$ containing $1$ has size $m\in\bn$.  Then, as in the proof of Theorem~\ref{eKn}, the restriction $\al|_A$ belongs to $C(\K_A)$ and the restriction $\al|_{\bn\sm A}$ belongs to $E(\K_{\bn\sm A})$.  But, since $\rank(\al)=r$, it follows that either
\bit
\item[(i)] $\al|_A\in C_0(\K_A)$ and $\al|_{\bn\sm A}\in E(D_r(\K_{\bn\sm A}))$, or
\item[(ii)] $\al|_A\in C_1(\K_A)$ and $\al|_{\bn\sm A}\in E(D_{r-1}(\K_{\bn\sm A}))$.
\eit
The proof concludes in a similar fashion to the proof of Theorem \ref{eKn}. \epf

As usual, for an odd integer~$k$, we write $k!!=k(k-2)\cdots3\cdot1$, and we interpret $(-1)!!=1$.

\bs
\begin{lemma}\label{lamKn}
If $\rho(\K_n)$ denotes the number of $\R$-classes in $D_0(\K_n)$ where $\K_n$ is one of $\P_n,\B_n,\PB_n$, then
\[
\rho(\P_n) = B(n) , \qquad\quad
\rho(\B_n) = \begin{cases}
(n-1)!! &\text{if $n$ is even}\\
0 &\text{if $n$ is odd,}
\end{cases} \qquad\quad
\rho(\PB_n) = a_n,
\]
where $B(n)$ is the $n$th Bell number, and $a_n$ satisfies the recurrence
\[
a_0=a_1=1 \COMMA a_n=a_{n-1}+(n-1)a_{n-2} \ \ \text{for $n\geq2$.}
\]
\end{lemma}

\pf The results concerning $\P_n$ and $\B_n$ are well-known; see for example the proof of \cite[Theorems~7.5 and~8.4]{EastGray}.  For the $\PB_n$ statement, note that, since $\dom(\al)=\emptyset$ for all $\al\in D_0(\PB_n)$, Theorem \ref{green} says that the $\R$- classes of $\PB_n$ are indexed by the equivalence relations $\ve$ on $\bn$ that satisfy the condition that each equivalence class has size $1$ or $2$.  Such equivalences are in one-one correspondence with the involutions (i.e., self-inverse permutations) of $\bn$, of which there are $a_n$ (see A000085 on the OEIS \cite{OEIS}). \epf

\begin{rem}
The \emph{completely regular} elements of a semigroup are those that are $\H$-related to an idempotent.  Because the $\H$-class of any idempotent from a (non-empty) $\D$-class $\DrKn$ is isomorphic to the symmetric group $\S_r$, it follows that the number of completely regular elements in $\K_n$ is equal to $\sum_{r=0}^n r!e(\DrKn)$.
\end{rem}

\subsection{The partition monoid}\label{subsect:Pn}

In this section, we obtain formulae for $c_0(\P_n),c_1(\P_n),c(\P_n)$.  Together with Theorems \ref{eKn}, \ref{eDrKn} and~\ref{eDrKn_rec}, this yields formulae and recurrences for $e(\P_n)$ and $e(\DrPn)$.  The key step is to enumerate the pairs of equivalence relations on $\bn$ with specified numbers of equivalence classes and whose join is equal to the universal relation $\bn\times\bn$.

Let $E(n)$ denote the set of all equivalence relations on $\bn$.  If $\ve\in E(n)$, we denote by $\bn/\ve$ the quotient of $\bn$ by $\ve$, which consists of all $\ve$-classes of $\bn$.  For $r,s\in\bn$, we define sets
\begin{align*}
E(n,r) &= \set{\ve\in E(n)}{|\bn/\ve|=r}, \\
E(n,r,s) &= \set{(\ve,\eta)\in E(n,r)\times E(n,s)}{\ve\vee\eta=\bn\times\bn},
\end{align*}
and we write $e(n,r,s)=|E(n,r,s)|$.

\bs
\begin{prop}\label{cPn}
If $n\geq1$, then 
\[
c_0(\P_n) = \sum_{r,s\in\bn} e(n,r,s) \COMMA
c_1(\P_n) = \sum_{r,s\in\bn} rs\cdot e(n,r,s) \COMMA
c(\P_n) = \sum_{r,s\in\bn} (1+rs)e(n,r,s).
\]
A recurrence for the numbers $e(n,r,s)$ is given in Proposition \ref{enrs}.
\end{prop}

\pf Let $r,s\in\bn$ and consider a pair $(\ve,\eta)\in E(n,r,s)$.  We count the number of idempotent partitions $\al\in C(\P_n)$ such that $\keru(\al)=\ve$ and $\kerl(\al)=\eta$.  Clearly there is a unique such $\al$ satisfying $\rank(\al)=0$.  To specify such an $\al$ with $\rank(\al)=1$, we must also specify one of the $\ve$-classes and one of the $\eta$-classes to form the unique transversal block of $\al$, so there are $rs$ of these.  Since there are $e(n,r,s)$ choices for $(\ve,\eta)$, the statements follow after summing over all~$r,s$.~\epf

\begin{rem}
The numbers $c_0(\P_n)$ count the number of pairs of equivalences on $\bn$ whose join is $\bn\times\bn$.  These numbers may be found in Sequence A060639 on the OEIS \cite{OEIS}.
\end{rem}

For the proof of the following result, we denote by $\ve_{ij}\in E(n)$ the equivalence relation whose only non-trivial equivalence class is $\{i,j\}$.  On a few occasions in the proof, we will make use of the (trivial) fact that if $\ve\in E(n,r)$, then $\ve\vee\ve_{ij}$ has at least $r-1$ equivalence classes.  As usual, we write $S(n,r)=|E(n,r)|$; these are the (unsigned) Stirling numbers of the second kind.

\bs
\begin{prop}\label{enrs}
The numbers $e(n,r,s)=|E(n,r,s)|$ satisfy the recurrence:
\begin{align*}
e(n,r,1) &= S(n,r)  \\
e(n,1,s) &= S(n,s)  \\
e(n,r,s) &= s\cdot e(n-1,r-1,s) + r\cdot e(n-1,r,s-1) + rs\cdot e(n-1,r,s) \\
& \qquad +\sum_{m=1}^{n-2} {n-2 \choose m} \sum_{a=1}^{r-1}\sum_{b=1}^{s-1} \big( a(s-b)+b(r-a)\big) e(m,a,b) e(n-m-1,r-a,s-b).\\
& &&\hspace{-5cm}\text{if $r,s\geq2$.}
\end{align*}
\end{prop}

\pf The $r=1$ and $s=1$ cases are clear, so suppose $r,s\geq2$.  Consider a pair $(\ve,\eta)\in E(n,r,s)$.  We consider several cases.  Throughout the proof, we will write $\bnf=\{1,\ldots,n-1\}$.

{\bf Case 1.}  Suppose first that $\{n\}$ is an $\ve$-class.  Let $\ve'=\ve\cap(\bnf\times\bnf)$ and $\eta'=\eta\cap(\bnf\times\bnf)$ denote the induced equivalence relations on $\bnf$.  Then we clearly have $\ve'\in E(n-1,r-1)$.  Also,~$\{n\}$ cannot be an $\eta$-class, or else then $\{n\}$ would be an $\ve\vee\eta$-class, contradicting the fact that $\ve\vee\eta=\bn\times\bn$.  It follows that $\eta'\in E(n-1,s)$.

Next we claim that $\ve'\vee\eta'=\bnf\times\bnf$.  Indeed, suppose to the contrary that $\ve'\vee\eta'$ has $k\geq2$ equivalence classes.  Let $\eta''\in E(n)$ be the equivalence on $\bn$ obtained from $\eta'$ by declaring $\{n\}$ to be an $\eta''$-class.  Then $\ve\vee\eta''$ has $k+1$ equivalence classes.  But $\eta=\eta''\vee\ve_{in}$ for some $i\in\bnf$.  It follows that $\ve\vee\eta=(\ve\vee\eta'')\vee\ve_{in}$ has $(k+1)-1=k\geq2$ equivalence classes, contradicting the fact that $\ve\vee\eta=\bn\times\bn$.  So this establishes the claim.

It follows that $(\ve',\eta')\in E(n-1,r-1,s)$.  So there are $e(n-1,r-1,s)$ such pairs.  We then have to choose which block of $\eta'$ to put $n$ into when creating $\eta$, and this can be done in $s$ ways.  So it follows that there are $s\cdot e(n-1,r-1,s)$ pairs $(\ve,\eta)$ in Case 1.

{\bf Case 2.}  By symmetry, there are $r\cdot e(n-1,r,s-1)$ pairs $(\ve,\eta)$ in the case that $\{n\}$ is an $\eta$-class.

{\bf Case 3.}  Now suppose that $\{n\}$ is neither an $\ve$-class nor an $\eta$-class.  Again, let $\ve',\eta'$ be the induced equivalences on $\bnf$.  This time, $\ve'\in E(n-1,r)$ and $\eta'\in E(n-1,s)$.  We now consider two subcases.

{\bf Case 3.1.}  If $\ve'\vee\eta'=\bnf\times\bnf$, then $(\ve',\eta')\in E(n-1,r,s)$.  By similar reasoning to that above, there are $rs\cdot e(n-1,r,s)$ pairs $(\ve,\eta)$ in this case.

{\bf Case 3.2.}  Finally, suppose $\ve'\vee\eta'\not=\bnf\times\bnf$, and denote by $k$ the number of $\ve'\vee\eta'$-classes.  We claim that $k=2$.  Indeed, suppose this is not the case.  By assumption, $k\not=1$, so it follows that $k\geq3$.  Let $\ve''$ and $\eta''$ be the equivalence relations on $\bn$ obtained from $\ve'$ and $\eta'$ by declaring $\{n\}$ to be an $\ve''$- and $\eta''$-class.  Then $\ve''\vee\eta''$ has $k+1$ equivalence classes, and $\ve=\ve''\vee\ve_{in}$ and $\eta=\eta''\vee\ve_{jn}$ for some $i,j\in\bnf$.  So
$
\ve\vee\eta = (\ve''\vee\eta'')\vee\ve_{in}\vee\ve_{jn}
$
has at least $(k+1)-2=k-1\geq2$ equivalence classes, a contradiction.  So we have proved the claim.

Denote by $B_1$ the $\ve'\vee\eta'$-class of $\bnf$ containing $1$, and let the other $\ve'\vee\eta'$-class be $B_2$, noting that $1\leq |B_1|\leq n-2$.  If $|B_1|=m$, then there are ${n-2\choose m-1}$ ways to choose $B_1$ (and $B_2=\bnf\sm B_1$ is then fixed).

For $i=1,2$, let $\ve_i=\ve\cap(B_i\times B_i)$ and $\eta_i=\eta\cap(B_i\times B_i)$.  Note that $\ve_i\vee\eta_i=B_i\times B_i$ for each $i$.  Let $a=|B_1/\ve_1|$ and $b=|B_1/\eta_1|$.  So $1\leq a\leq r-1$ and $1\leq b\leq s-1$, and also $|B_2/\ve_2|=r-a$ and $|B_2/\eta_2|=s-b$.  So, allowing ourselves to abuse notation slightly, we have $(\ve_1,\eta_1)\in E(m,a,b)$ and $(\ve_2,\eta_2)\in E(n-m-1,r-a,s-b)$.  So there are $e(m,a,b) e(n-m-1,r-a,s-b)$ ways to choose $\ve_1,\ve_2,\eta_1,\eta_2$.

We must also choose which blocks of $\ve'$ and $\eta'$ to add $n$ to, when creating $\ve,\eta$ from $\ve',\eta'$.  But, in order to ensure that $\ve\vee\eta=\bn\times\bn$, if we add $n$ to one of the $\ve'$-classes in $B_1$, we must add $n$ to one of the $\eta'$-classes in $B_2$, and vice versa.  So there are $a(s-b)+b(r-a)$ choices for the blocks to add $n$ to.

Multiplying the quantities obtained in the previous three paragraphs, and summing over the appropriate values of $m,a,b$, we get a total of
\[
\sum_{m=1}^{n-2} {n-2 \choose m} \sum_{a=1}^{r-1}\sum_{b=1}^{s-1} \big( a(s-b)+b(r-a)\big) e(m,a,b) e(n-m-1,r-a,s-b)
\]
pairs $(\ve,\eta)$ in Case 3.2.

Adding the values from all the above cases gives the desired result. \epf

\bsa
\subsection{The Brauer monoid}\label{subsect:Bn}

We now apply the general results above to derive a formula for $e(\B_n)$.  As in the previous section, the key step is to obtain formulae for $c_0(\B_n),c_1(\B_n),c(\B_n)$, but the simple form of these values (see Proposition \ref{cBn}) allows us to obtain neat expressions for $e(\B_n)$ and $e(\DrBn)$ (see Theorem~\ref{eBn}).
But first it will be convenient to prove a result concerning the graphs $\Ga(\al)$, where $\al$ belongs to the larger partial Brauer monoid $\PB_n$, as it will be useful on several occasions (these graphs were defined after Theorem \ref{green}).

\bs
\begin{lemma}\label{Gaalform}
Let $\al\in C(\PB_X)$ where $X$ is finite.  Then $\Ga(\al)$ is either a cycle or a path.
\end{lemma}

\pf The result is trivial if $|X|=1$, so suppose $|X|\geq2$.  In the graph $\Ga(\al)$, no vertex can have two red or two blue edges coming out of it, so it follows that the degree of each vertex is at most~$2$.
It follows that $\Ga(\al)$ is a union of paths and cycles.  Since $\Ga(\al)$ is connected, we are done. \epf

\begin{prop}\label{cBn}
If $n\geq1$, then
\begin{align*}
c_0(\B_n) &= \begin{cases}
0 &\text{if $n$ is odd}\\
(n-1)! &\text{if $n$ is even,}
\end{cases}\\
c_1(\B_n)&= \begin{cases}
n! &\hspace{0.95cm}\text{if $n$ is odd}\\
0 &\hspace{0.95cm}\text{if $n$ is even,}
\end{cases}\\
c(\B_n) &= \begin{cases}
n! &\text{if $n$ is odd}\\
(n-1)! &\text{if $n$ is even.}
\end{cases}
\end{align*}
\end{prop}

\pf Let $\al$ be an irreducible idempotent from $\B_n$.  By Lemma \ref{tech3}, and the fact that $\rank(\be)\in\{n,n-2,n-4,\ldots\}$ for all $\be\in\B_n$, we see that $\rank(\al)=0$ if $n$ is even, while $\rank(\al)=1$ if $n$ is odd.
If $n$ is even, then by Lemma \ref{Gaalform} and the fact that $\Ga(\al)$ has the same number of red and blue edges, $\Ga(\al)$ is a cycle $1{\red-}i_2{\blue-}i_3{\red-}\cdots{\red-}i_n{\blue-}1$, where $\{i_2,\ldots,i_n\}=\{2,\ldots,n\}$, and there are precisely $(n-1)!$ such cycles.  Similarly, if $n$ is odd, then $\Ga(\al)$ is a path $i_1{\red-}i_2{\blue-}i_3{\red-}\cdots{\blue-}i_{n}$, where $\{i_1,\ldots,i_n\}=\bn$, and there are $n!$ such paths. \epf

\begin{thm}\label{eBn}
Let $n\in\mathbb N$ and put $k=\lfloor \frac n2\rfloor$.  Then
\begin{align*}
e(\B_n) &= \sum_{\mu\vdash n} \frac{n!}{\prod_{i=1}^n\mu_i! \cdot \prod_{j=1}^{k}(2j)^{\mu_{2j}}}.
\intertext{If $0\leq r\leq n$, then}
e(\DrBn) &= \sum_{\mu\vdash n \atop \mu_1+\mu_3+\cdots=r} \frac{n!}{\prod_{i=1}^n\mu_i! \cdot \prod_{j=1}^{k}(2j)^{\mu_{2j}}}.
\end{align*}
\end{thm}

\pf By Theorem \ref{eKn} and Proposition \ref{cBn},
\begin{align*}
e(\B_n) = n!\cdot \sum_{\mu\vdash n} \frac{(1!)^{\mu_1}(1!)^{\mu_2}(3!)^{\mu_3}(3!)^{\mu_4}\cdots}{\mu_1!\cdots\mu_n! \cdot (1!)^{\mu_1}(2!)^{\mu_2}(3!)^{\mu_3} (4!)^{\mu_4}\cdots}
= n!\cdot \sum_{\mu\vdash n} \frac{1}{\mu_1!\cdots\mu_n! \cdot 2^{\mu_2}\cdot 4^{\mu_4}\cdots},
\end{align*}
establishing the first statement.  For the second, suppose $\mu=(m_1,\ldots,m_k)=(1^{\mu_1},\ldots,n^{\mu_n})\vdash n$.  Theorem \ref{eDrKn} gives
\[
e_\mu(\DrBn) = 
 \frac{n!}{\prod_{i=1}^n\mu_i!(i!)^{\mu_i}}
 \sum_{A\sub\bk \atop |A|=r} \left(\prod_{i\in A}c_1(\B_{m_i})\cdot\prod_{j\in\bk\sm A}c_0(\B_{m_j})\right).
\]
By Proposition \ref{cBn}, we see that for $A\sub\bk$ with $|A|=r$, 
\begin{align*}
\text{$\prod_{i\in A}c_1(\B_{m_i})\cdot\prod_{j\in\bk\sm A}c_0(\B_{m_j})\not=0$} &\iff \text{$m_i$ is odd for all $i\in A$ and $m_j$ is even for all $j\in\bk\sm A$} \\
&\iff \text{$A=\set{i\in\bk}{\text{$m_i$ is odd}}$.}
\intertext{So}
\text{$e_\mu(\DrBn)\not=0$} &\iff \text{$\set{i\in\bk}{\text{$m_i$ is odd}}$ has size $r$} \\
&\iff \text{$\mu_1+\mu_3+\cdots=r$,}
\end{align*}
in which case,
\[
e_\mu(\DrBn) = n!\cdot  \frac{(1!)^{\mu_1}(1!)^{\mu_2}(3!)^{\mu_3}(3!)^{\mu_4}\cdots}{\mu_1!\cdots\mu_n! \cdot (1!)^{\mu_1}(2!)^{\mu_2}(3!)^{\mu_3} (4!)^{\mu_4}\cdots}
= n!\cdot \frac{1}{\mu_1!\cdots\mu_n! \cdot 2^{\mu_2}\cdot 4^{\mu_4}\cdots}.
\]
Summing over all $\mu\vdash n$ with $\mu_1+\mu_3+\cdots=r$ gives the desired expression for $e(\DrBn)$. \epf

\begin{rem}
The formula for $e(\B_n)$ may be deduced from \cite[Proposition 4.10]{Larsson}.  Note that $e(\DrBn)\not=0$ if and only if $n\equiv r\pmod2$.
\end{rem}

Proposition \ref{cBn} also leads to a simple form of the recurrences from Theorems \ref{eKn} and \ref{eDrKn_rec} for the numbers $e(\B_n)$ and $e(\DrBn)$.

\bs
\begin{thm}
The numbers $e(\B_n)$ satisfy the recurrence:
\begin{align*}
e(\B_0)=1 \COMMA e(\B_n) &= \sum_{i=1}^{\lfloor \frac n2\rfloor} {n-1\choose 2i-1}(2i-1)! \;\:\!\! e(\B_{n-2i})\\
\epfreseq
& \qquad\qquad\qquad + \sum_{i=0}^{\lfloor \frac {n-1}2\rfloor} {n-1\choose 2i}(2i+1)! \;\:\!\! e(\B_{n-2i-1}) \qquad\text{for $n\geq1$.}
\end{align*}
\end{thm}

\bs
\begin{thm}
The numbers $e(\DrBn)$ satisfy the recurrence:
\begin{align*}
e(D_n(\B_n)) &= 1 \\ 
e(D_0(\B_n)) &= \begin{cases}
(n-1)!!^2 &\text{if $n$ is even}\\
0 &\text{if $n$ is odd}
\end{cases}  \\ 
e(\DrBn) &= \sum_{i=1}^{\lfloor \frac n2\rfloor} {n-1\choose 2i-1}(2i-1)! \;\:\!\! e(D_r(\B_{n-2i}))
\\ \epfreseq & \qquad\qquad\qquad
+ \sum_{i=0}^{\lfloor \frac {n-1}2\rfloor} {n-1\choose 2i}(2i+1)! \;\:\!\! e(D_{r-1}(\B_{n-2i-1})) \qquad\text{for $n\geq1$.}
\end{align*}
\end{thm}

\bsa
\subsection{The partial Brauer monoid}\label{subsect:PBn}

As usual, the key step in calculating $e(\PB_n)$ is to obtain formulae for $c(\PB_n)$.

\bs
\begin{prop}\label{cPBn}
If $n\geq1$, then
\begin{align*}
c_0(\PB_n) &= \begin{cases}
n! &\text{if $n$ is odd}\\
(n+1)\cdot(n-1)! &\text{if $n$ is even,}
\end{cases}\\
c_1(\PB_n) &= \begin{cases}
n! &\hspace{2.43 cm}\text{if $n$ is odd}\\
0 &\hspace{2.43 cm}\text{if $n$ is even,}
\end{cases}\\
c(\PB_n)&=\begin{cases}
2\cdot n! &\text{if $n$ is odd}\\
(n+1)\cdot (n-1)! &\text{if $n$ is even.}
\end{cases}
\end{align*}
\end{prop}

\pf Let $\al$ be an irreducible idempotent from $\PB_n$.  Suppose first that $n$ is odd.  By Lemma~\ref{Gaalform}, whether $\rank(\al)$ is equal to $0$ or $1$, $\Ga(\al)$ is a path $i_1{\red-}i_2{\blue-}i_3{\red-}\cdots{\blue-}i_{n}$, and there are $n!$ such paths.
Now suppose $n$ is even.  Then $\Ga(\al)$ is either a cycle $1{\red-}i_2{\blue-}i_3{\red-}\cdots{\red-}i_n{\blue-}1$, of which there are $(n-1)!$, or else a path $i_1{\red-}i_2{\blue-}i_3{\red-}\cdots{\red-}i_{n}$ or $i_1{\blue-}i_2{\red-}i_3{\blue-}\cdots{\blue-}i_{n}$, of which there are $n!/2$ of both kinds.  All of these have $\rank(\al)=0$, and adding them gives $n!+(n-1)!=(n+1)\cdot(n-1)!$.~\epf

\nc{\Bmu}{B_\mu}

\begin{thm}\label{ePBn}
Let $n\in\mathbb N$ and put $k=\lfloor \frac n2\rfloor$.  Then
\begin{align*}
e(\PB_n) &= n! \cdot\sum_{\mu\vdash n} \frac{\prod_{j=1}^k(1+\frac1{2j})^{\mu_{2j}}}{\prod_{i=1}^n\mu_i!} 
2^{\mu_1+\mu_3+\cdots}.
\intertext{If $0\leq r\leq n$, then}
e(\DrPBn) &= n! \cdot \!\!\!\!\!\!\! \sum_{\mu\vdash n \atop \mu_1+\mu_3+\cdots\geq r} \!\!\!\!\!\!\! \frac{\prod_{j=1}^k(1+\frac1{2j})^{\mu_{2j}}}{\prod_{i=1}^n\mu_i!}
{\mu_1+\mu_3+\cdots \choose r}.
\end{align*}
\end{thm}

\pf By Theorem \ref{eKn} and Proposition \ref{cPBn},
\begin{align*}
e(\PB_n) &= n!\cdot \sum_{\mu\vdash n} \frac{(2\cdot1!)^{\mu_1}(2\cdot3!)^{\mu_3}\cdots(3\cdot1!)^{\mu_2}(5\cdot3!)^{\mu_4}\cdots}{\mu_1!\cdots\mu_n! \cdot (1!)^{\mu_1}(3!)^{\mu_3} \cdots (2!)^{\mu_2}(4!)^{\mu_4}\cdots} \\
&= n!\cdot \sum_{\mu\vdash n} \frac{2^{\mu_1+\mu_3+\cdots}}{\mu_1!\cdots\mu_n! }
\left(\frac32\right)^{\mu_2}\left(\frac54\right)^{\mu_4} \cdots,
\end{align*}
giving the first statement.  For the second, suppose $\mu=(m_1,\ldots,m_k)=(1^{\mu_1},\ldots,n^{\mu_n})\vdash n$.  Theorem~\ref{eDrKn} gives
\[
e_\mu(\DrPBn) = 
 \frac{n!}{\prod_{i=1}^n\mu_i!(i!)^{\mu_i}}
 \sum_{A\sub\bk \atop |A|=r} \left(\prod_{i\in A}c_1(\PB_{m_i})\cdot\prod_{j\in\bk\sm A}c_0(\PB_{m_j})\right).
\]
Let $\Bmu=\set{i\in\bk}{\text{$m_i$ is odd}}$.  By Proposition \ref{cPBn}, we see that for $A\sub\bk$ with $|A|=r$, 
\[
\text{$\prod_{i\in A}c_1(\PB_{m_i})\cdot\prod_{j\in\bk\sm A}c_0(\PB_{m_j})\not=0$} \iff \text{$m_i$ is odd for all $i\in A$}
\iff \text{$A\sub\Bmu$.}
\]
In particular, $e_\mu(\DrPBn)\not=0$ if and only if $\mu_1+\mu_3+\cdots=|B_\mu|\geq r$.  For such a $\mu\vdash n$ and for $A\sub\Bmu$ with $|A|=r$,
\begin{align*}
\prod_{i\in A}c_1(\PB_{m_i})\cdot\prod_{j\in\bk\sm A}c_0(\PB_{m_j}) &= \prod_{i\in A}c_1(\PB_{m_i})\cdot\prod_{i\in \Bmu\sm A}c_0(\PB_{m_i})\cdot\prod_{j\in\bk\sm \Bmu}c_0(\PB_{m_j}) \\
&= \prod_{i\in \Bmu}m_i!\cdot\prod_{j\in\bk\sm \Bmu}(m_j+1)\cdot(m_j-1)! \\
&= (1!)^{\mu_1}(3!)^{\mu_3}\cdots(3\cdot 1!)^{\mu_2}(5\cdot3!)^{\mu_4}\cdots.
\end{align*}
Since there are ${\mu_1+\mu_3+\cdots \choose r}$ subsets $A\sub\Bmu$ with $|A|=r$, it follows that
\begin{align*}
e_\mu(\DrPBn) &= {\mu_1+\mu_3+\cdots \choose r} \cdot n! \cdot \frac{(1!)^{\mu_1}(3!)^{\mu_3}\cdots(3\cdot 1!)^{\mu_2}(5\cdot3!)^{\mu_4}\cdots}{\mu_1!\cdots\mu_n! \cdot (1!)^{\mu_1}(3!)^{\mu_3}\cdots(2!)^{\mu_2}(4!)^{\mu_4}\cdots} \\
&= {\mu_1+\mu_3+\cdots \choose r}
\frac{n!}{\mu_1!\cdots\mu_n!}
\left(\frac32\right)^{\mu_2}\left(\frac54\right)^{\mu_4} \cdots.
\end{align*}
Summing over all $\mu\vdash n$ with $\mu_1+\mu_3+\cdots\geq r$ gives the required expression for $e(\DrPBn)$. \epf

Again, the recurrences for the numbers $e(\PB_n)$ and $e(\DrPBn)$ given by Theorems \ref{eKn} and \ref{eDrKn_rec} take on a neat form.

\bs
\begin{thm}
The numbers $e(\PB_n)$ satisfy the recurrence:
\begin{align*}
e(\PB_0)=1 \COMMA e(\PB_n) &= \sum_{i=1}^{\lfloor \frac n2\rfloor} {n-1\choose 2i-1}(2i+1)\cdot(2i-1)! \;\:\!\! e(\PB_{n-2i}) \\
\epfreseq
& \qquad\qquad\qquad + 2\cdot\sum_{i=0}^{\lfloor \frac {n-1}2\rfloor} {n-1\choose 2i}(2i+1)! \;\:\!\! e(\PB_{n-2i-1}) \qquad\text{for $n\geq1$.}
\end{align*}
\end{thm}

\bs
\begin{thm}
The numbers $e(\DrPBn)$ satisfy the recurrence:
\begin{align*}
e(D_n(\PB_n)) &= 1 \\ 
e(D_0(\PB_n)) &= \begin{cases}
a_n^2 &\text{if $n$ is even}\\
0 &\text{if $n$ is odd}
\end{cases}  \\ 
e(\DrPBn) &= \sum_{i=1}^{\lfloor \frac n2\rfloor} {n-1\choose 2i-1}(2i+1)\cdot(2i-1)! \;\:\!\! e(D_r(\PB_{n-2i}))
\\
& \qquad
+ \sum_{i=0}^{\lfloor \frac {n-1}2\rfloor} {n-1\choose 2i}(2i+1)!\Big(e(D_{r}(\PB_{n-2i-1}))+e(D_{r-1}(\PB_{n-2i-1}))\Big) \quad\text{for $n\geq1$,}
\end{align*}
where the numbers $a_n$ are defined in Lemma \ref{lamKn}. \epfres
\end{thm}

\bsa
\subsection{Other subsemigroups}

We conclude this section with the observation that the general results above (Theorems~\ref{eKn},~\ref{eDrKn},~\ref{eDrKn_rec}) apply to many other subsemigroups of $\P_n$ (though the initial conditions need to be slightly modified in Theorem \ref{eDrKn_rec}).  As observed in \cite{JEgrpm,JEpnsn,EF}, the full transformation semigroup and the symmetric and dual symmetric inverse semigroups $\T_n,\I_n,\I_n^*$ are all (isomorphic to) subsemigroups of $\P_n$:
\bit
\item $\T_n\cong\set{\al\in\P_n}{\text{$\dom(\al)=\bn$ and $\kerl(\al)=\Delta$}}$,
\item $\I_n\cong\set{\al\in\P_n}{\text{$\keru(\al)=\kerl(\al)=\Delta$}}$,
\item $\I_n^*\cong\set{\al\in\P_n}{\text{$\dom(\al)=\codom(\al)=\bn$}}$,
\eit
where $\Delta=\set{(i,i)}{i\in\bn}$ denotes the trivial equivalence (that is, the equality relation), and the above mentioned theorems apply to these subsemigroups.  For example, one may easily check that $c(\T_n)=c_1(\T_n)=n$, so that Theorem \ref{eKn} gives rise to the formula
$$
e(\T_n) = n!\cdot\sum_{\mu\vdash n} \prod_{i=1}^n \frac{1}{\mu_i! ((i-1)!)^{\mu_i}},
$$
and the recurrence
$$
e(\T_0)=1, \quad e(\T_n) = \sum_{m=1}^n {n-1\choose m-1} \cdot m\cdot e(\T_{n-m}) \quad\text{for $n\geq1$.}
$$
As noted in the Introduction, the usual formula is $e(\T_n)=\sum_{k=1}^n{n\choose k}k^{n-k}$.  The recurrence for $e(\I_n^*)$, combined with the fact that $e(\I_n^*)=B(n)$ is the $n$th Bell number \cite{FL}, leads to
\[
B(n+1) = \sum_{k=0}^n{n\choose k}B(k),
\]
a well-known identity.
We leave it to the reader to explore further if they wish.

\section{A different approach for $\B_n$ and $\PB_n$}\label{sect:diff}

We now outline an alternative method for determining $e(\B_n)$ and $e(\PB_n)$.  This approach will also allow us to determine the number of idempotents in an arbitrary $\R$-, $\L$- and $\D$-class of $\B_n$ and $\PB_n$.  One advantage of this method is that we do not need to take sums over integer partitions; rather, everything depends on sequences defined by some fairly simple recurrence relations (see Theorems \ref{anr} and \ref{anrt}).  The key idea is to define a variant of the graph $\Ga(\al)$ in the case of $\al\in\PB_n$.

Let $\al\in\PB_X$.  We define $\Lamu(\al)$ (resp.~$\Laml(\al)$) to be the graph obtained from $\Gau(\al)$ (resp.~$\Gal(\al)$) by adding a red (resp.~blue) loop at each vertex $i\in X$ if $\{i\}$ (resp.~$\{i'\}$) is a block of $\al$.  And we define $\Lam(\al)$ to be the graph with vertex set $X$ and all the edges from both $\Lamu(\al)$ and $\Laml(\al)$.  Some examples are given in Figure~\ref{fig:Lam} with $X$ finite.  Note that $\Lam(\al)=\Ga(\al)$ if and only if $\al\in\B_X$.

\begin{figure}[ht]
\begin{center}
\begin{tikzpicture}[xscale=.7,yscale=0.7]
    \uv0
    \uv1
    \uv2
    \uv3
    \uv4
    \uv5
    \lv0
    \lv1
    \lv2
    \lv3
    \lv4
    \lv5
    \arcup01
    \arcdn03
    \arcdn45
    \cve31
    \larcr01
    \larcb03
    \larcb45
         \llloopr2{.2}
         \lrloopb2{.2}
         \ldloopr4{.2}
         \ldloopr5{.2}
    \llv0
    \llv1
    \llv2
    \llv3
    \llv4
    \llv5
\end{tikzpicture}
\qquad\qquad
\begin{tikzpicture}[xscale=.7,yscale=0.7]
    \uv0
    \uv1
    \uv2
    \uv3
    \uv4
    \uv5
    \lv0
    \lv1
    \lv2
    \lv3
    \lv4
    \lv5
    \arcup01
    \arcdn02
    \arcdn45
    \cve31
    \larcr01
    \larcb02
    \larcb45
         \ldloopr2{.2}
         \luloopb3{.2}
         \ldloopr4{.2}
         \ldloopr5{.2}
    \llv0
    \llv1
    \llv2
    \llv3
    \llv4
    \llv5
\end{tikzpicture}
\qquad\qquad
\begin{tikzpicture}[xscale=.7,yscale=0.7]
    \uv0
    \uv1
    \uv2
    \uv3
    \uv4
    \uv5
    \lv0
    \lv1
    \lv2
    \lv3
    \lv4
    \lv5
    \arcup02
    \arcup45
    \arcdn23
    \arcdn45
    \cve10
    \cve31
    \larcr02
    \larcr45
    \larcb23
    \larcb45
    \llv0
    \llv1
    \llv2
    \llv3
    \llv4
    \llv5
\end{tikzpicture}
\end{center}
\caption{Elements $\al,\be,\ga$ (left to right) of the partial Brauer monoid $\PB_6$ and their graphs $\Lam(\al),\Lam(\be),\Lam(\ga)$ (below).}
\label{fig:Lam}
\end{figure}

Since the graph $\Lamu(\al)$ (resp.~$\Laml(\al)$) determines (and is determined by) $\dom(\al)$ and $\keru(\al)$ (resp.~$\codom(\al)$ and $\kerl(\al)$), we immediately obtain the following from Theorem \ref{green}.

\bs
\begin{cor}\label{greenPBn}
Let $X$ be any set (finite or infinite).
For each $\al,\be\in\PB_X$, we have
\ben
\item[\emph{(i)}] $\al\R\be$ if and only if $\Lamu(\al)=\Lamu(\be)$,
\item[\emph{(ii)}] $\al\L\be$ if and only if $\Laml(\al)=\Laml(\be)$,
\item[\emph{(iii)}] $\al\H\be$ if and only if $\Lam(\al)=\Lam(\be)$. \epfres
\een
\end{cor}

We now aim to classify the graphs on vertex set $X$ that are of the form $\Lam(\al)$ for some $\al\in\PB_X$, and we will begin with the irreducible idempotents.

\bs
\begin{lemma}\label{Lamalform}
Let $\al\in C(\PB_X)$ where $X$ is finite.  Then $\Lam(\al)$ is of one of the following four forms:
\ben
\item[\emph{(1)}] $\textlinegraph{.3}{
\bluto12
\redto23
\bluto3{3.5}
\redto{4.5}5
\dotto{3.5}{4.5}
\overt1
\overt2
\overt3
\overt5
}$ \emph{:} an alternating-colour path of even length,
\item[\emph{(2)}] $\textlinegraph{.3}{
\bluto12
\redto23
\bluto3{3.5}
\bluto{4.5}5
\dotto{3.5}{4.5}
\arcr15
\overt1
\overt2
\overt3
\overt5
}$ \emph{:} an alternating-colour circuit of even length,
\item[\emph{(3)}] $\textlinegraph{-.2}{
\bluto12
\redto23
\bluto3{3.5}
\redto{4.5}5
\dotto{3.5}{4.5}
\lloopr1{.2}
\rloopb5{.2}
\overt1
\overt2
\overt3
\overt5
}$ \emph{:} an alternating-colour path of even length with loops,
\item[\emph{(4)}] $\textlinegraph{-.2}{
\bluto12
\redto23
\bluto3{3.5}
\bluto{4.5}5
\dotto{3.5}{4.5}
\lloopr1{.2}
\rloopr5{.2}
\overt1
\overt2
\overt3
\overt5
}$ \ or \ $\textlinegraph{-.2}{
\redto12
\bluto23
\redto3{3.5}
\redto{4.5}5
\dotto{3.5}{4.5}
\lloopb1{.2}
\rloopb5{.2}
\overt1
\overt2
\overt3
\overt5
}$ \emph{:} an alternating-colour path of odd length with loops.
\een
If $\al\in C(\B_X)$, then $\Lam(\al)$ is of the form \emph{(1)} or \emph{(2)}.
\end{lemma}

\pf By Lemma \ref{Gaalform}, we know that $\Ga(\al)$ is either a cycle or a path.  If $\Ga(\al)$ is a cycle, then $\Lam(\al)=\Ga(\al)$ is of type~(2).  Next suppose $\Ga(\al)$ is a path, and write $n=|X|$.
We consider the case in which $n$ is odd (so the path is of even length).  Re-labelling the elements of $X$ if necessary, we may assume that
$\Ga(\al)$ is the path $1{\blue-}2{\red-}3{\blue-}\cdots{\red-}n$.
Since $\Ga(\al)$ completely determines $\keru(\al)$ and $\kerl(\al)$, it follows that $\al$ must be one of
$\disppartnthreeone{\uarc23 \uarc56 \darc12 \darc45 \cvs16 }$ or $\disppartnthreeone{\uarc23 \uarc56 \darc12 \darc45}$, in which case $\Lam(\al)$ is of type (1) or (3), respectively.  (Note that the preceeding discussion include the case $n=1$, where we have $\Ga(\al)=\textlinegraph{.3}{\overt1}$, so the two possibilities for $\Lam(\al)$ are $\textlinegraph{.3}{\overt1}$ or $\textlinegraph{-.2}{\lloopr1{.2}\rloopb1{.2}\overt1}$, which are of type (1) and (3), respectively.)  The case in which $n$ is even is similar, and leads to $\Lam(\al)$ being of type~(4).  The statement concerning $\B_X$ is clear, seeing as elements of $\B_X$ have no singleton blocks.~\epf

Now consider a graph $\Lam$ with edges coloured red or blue.  We say that $\Lam$ is \emph{balanced} if it is a disjoint union of finitely many subgraphs of types (1--4) from Lemma~\ref{Lamalform}.  We call a balanced graph $\Lam$ \emph{reduced} if it is a disjoint union of finitely many subgraphs of types (1--2) from Lemma~\ref{Lamalform}. If $X$ is a finite set, we write $\Bal(X)$ (resp.~$\Red(X)$) for the set of all balanced (resp.~reduced balanced) graphs with vertex set $X$.

\bs
\begin{prop}\label{bijections}
If $X$ is a finite set, then the maps
\[
\Phi:E(\PB_X)\to\Bal(X):\al\mt\Lam(\al) \quad\text{and}\quad \Psi=\Phi|_{E(\B_X)}:E(\B_X)\to\Red(X):\al\mt\Lam(\al)=\Ga(\al)
 \]
are bijections.  If $\al\in E(\PB_X)$, then $\rank(\al)$ is equal to the number of connected components of $\Lam(\al)$ of type \emph{(1)} as listed in Lemma \ref{Lamalform}.
\end{prop}

\pf Let $\al\in E(\PB_X)$, and write $\al=\al_1\oplus\cdots\oplus\al_k$ where $\al_1,\ldots,\al_k$ are the irreducible components of $\al$.  Then $\Lam(\al)$ is the disjoint union of the subgraphs $\Lam(\al_1),\ldots,\Lam(\al_k)$, and is therefore reduced, by Lemma \ref{Lamalform}.  If, in fact, $\al\in E(\B_X)$, then each of $\Lam(\al_1),\ldots,\Lam(\al_k)$ must be of the form (1) or~(2), since $\al$ has no blocks of size $1$.  This shows that $\Phi$ and $\Psi$ do indeed map $E(\PB_X)$ and $E(\B_X)$ to $\Bal(X)$ and $\Red(X)$, respectively.  Note also that $\rank(\al)=\rank(\al_1)+\cdots+\rank(\al_k)$ is equal to the number of rank $1$ partitions among $\al_1,\ldots,\al_k$, and that the rank of some $\be\in C(\PB_Y)$ is equal to $1$ if and only if $\Lam(\be)$ is of type (1).

Let $\Lam\in\Bal(X)$, and suppose $\Lam_1,\ldots,\Lam_k$ are the connected components of $\Lam$, with vertex sets $X_1,\ldots,X_k$, respectively.  Then there exist irreducible idempotents $\al_i\in C(\P_{X_i})$ with ${\Lam(\al_i)=\Lam_i}$ for each $i$, and it follows that $\Lam=\Lam(\al_1\oplus\cdots\oplus\al_k)$, showing that $\Phi$ is surjective.  If $\Lam\in\Red(X)$, then $\al_1\oplus\cdots\oplus\al_k\in\B_X$.
Finally, if $\al,\be\in E(\PB_X)$ are such that $\Lam(\al)=\Lam(\be)$, then $\al\H\be$ by Corollary \ref{greenPBn}, so that $\al=\be$ (as $\H$ is idempotent-separating), whence $\Phi$ (and hence also $\Psi$) is injective. \epf

\bsa
\subsection{The Brauer monoid}

For $\al\in\B_n$, we write
\[
\RaBn=\set{\be\in\B_n}{\Lamu(\be)=\Lamu(\al)} \AND \LaBn=\set{\be\in\B_n}{\Laml(\be)=\Laml(\al)}.
\]
By Corollary \ref{greenPBn}, these are precisely the $\R$- and $\L$-classes of $\al$ in $\B_n$.
At this point, it will be convenient to introduce an indexing set.  Put $\bn^0=\bn\cup\{0\}$, and let $I(n)=\set{r\in\bn^0}{n-r\in2\mathbb Z}$.  For $r\in I(n)$, let
\[
\DrBn = \set{\al\in\B_n}{\rank(\al)=r}.
\]
By Theorem \ref{green}, we see that these are precisely the $\D$-classes of $\B_n$.
We will need to know the number of $\R$-classes (which is equal to the number of $\L$-classes) in a given $\D$-class of $\B_n$.

\begin{lemma}[{See the proof of \cite[Theorem 8.4]{EastGray}}]\label{rhonr}
For $n\in\N$ and $r=n-2k\in I(n)$, the number of $\R$-classes (and $\L$-classes) in the $\D$-class $\DrBn$ is equal to
\[
\epfreseq
\rho_{nr}={n\choose r}(2k-1)!!=\frac{n!}{2^k k! r!}.
\]
\end{lemma}

\bs
\begin{thm}\label{anr}
Define a sequence $a_{nr}$, for $n\in\mathbb N$ and $r\in I(n)$, by
\begin{align*}
a_{nn} &= 1 &&\text{for all $n$}\\
a_{n0} &= (n-1)!! &&\text{if $n$ is even}\\
a_{nr} &= a_{n-1,r-1}+(n-r)a_{n-2,r} &&\text{if $1\leq r\leq n-2$.}
\end{align*}
Then for any $n\in\N$ and $r\in I(n)$, and with $\rho_{nr}$ as in Lemma \ref{rhonr}:
\ben
\item[\emph{(i)}] $e(\RaBn)=e(\LaBn)=a_{nr}$ for any $\al\in \DrBn$,
\item[\emph{(ii)}] $e(\DrBn)=\rho_{nr}a_{nr}$,
\item[\emph{(iii)}] $e(\B_n)=\sum_{r\in I(n)}\rho_{nr}a_{nr}$.
\een
\end{thm}

\pf Note that (iii) follows from (ii), which follows from (i) and Lemma \ref{rhonr}, so it suffices to prove~(i).  Let $\al\in \DrBn$.
Re-labelling the points from $\bn$, if necessary, we may assume that
\[
\Gau(\al) = \Lamu(\al) = \textlinegraph{-5.0}{
\redto45
\redto78
\dotto{1.4}{2.6}
\dotto{5.4}{6.6}
\overt1
\overt3
\overt4
\overt5
\overt7
\overt8
\draw(2,0)node[below]{$\underbrace{\phantom{mmmm...}}_r$};
}\ .
\vspace{-\bigskipamount} 
\]
Let $A_{nr}$ be the set of reduced balanced graphs on vertex set $\bn$ with the same red edges as $\Gau(\al)$.  By Corollary \ref{greenPBn} and Proposition \ref{bijections}, $|A_{nr}|=e(\RaBn)$.  Put $a_{nr}=|A_{nr}|$.  We show that $a_{nr}$ satisfies the stated recurrence.  By symmetry, $e(\LaBn)=a_{nr}$.

Clearly $a_{nn}=1$ for all $n$.  If $n$ is even, then $a_{n0}$ is the number of ways to match the vertices from $\bn$ with $n/2$ non-intersecting (blue) arcs, which is equal to $(n-1)!!$.  Suppose now that $1\leq r\leq n-2$.  Elements of $A_{nr}$ come in two kinds:
\ben 
\item those for which $1$ is a connected component of its own, and
\item those for which $1$ is an endpoint of an even length alternating path.
\een
There are clearly $a_{n-1,r-1}$ elements of $A_{nr}$ of type 1.  Suppose now that $\Ga\in A_{nr}$ is a graph of type 2.  There are $n-r$ possible vertices for vertex $1$ to be joined to by a blue edge.  Suppose the vertex adjacent to $1$ is $x$.  Removing these two vertices, as well as the blue edge $1{\blue-}x$ and the red edge adjacent to $x$ (and relabelling the remaining vertices), yields an element of $A_{n-2,r}$.  Since this process is reversible, there are $(n-r)a_{n-2,r}$ elements of $A_{nr}$ of type 2.  Adding these gives $a_{nr} = a_{n-1,r-1}+(n-r)a_{n-2,r}$. \epf

\bsa
\subsection{The partial Brauer monoid}

For $\al\in\PB_n$, we write
\[
\RaPBn=\set{\be\in\PB_n}{\Lamu(\be)=\Lamu(\al)} \AND \LaPBn=\set{\be\in\PB_n}{\Laml(\be)=\Laml(\al)}.
\]

By Corollary \ref{greenPBn}, these are precisely the $\R$- and $\L$-classes of $\al$ in $\PB_n$.  For $r\in\bn^0$, let
\[
\DrPBn =  \set{\al\in\PB_n}{\rank(\al)=r}.
\]
Again, these are precisely the $\D$-classes of $\PB_n$.  But unlike the case of $\B_n$, it is not true that any two $\D$-related elements of $\PB_n$ are $\R$-related to the same number of idempotents.  So we will obtain a formula for $e(\RaPBn)$ that will depend on the paramaters $r,t$, where $r=\rank(\al)$ and $t$ is the number of singleton non-transversal upper-kernel classes.  Note that $n,r,t$ are constrained by the requirement that $n-r-t$ is even.  With this in mind, we define an indexing set
\[
J(n)=\set{(r,t)\in\bn^0\times\bn^0}{t\in I(n-r)} = \set{(r,t)\in\bn^0\times\bn^0}{n-r-t\in2\mathbb Z}.
\]
There is a dual statement of the following lemma, but we will not state it.

\begin{lemma}\label{rhonrt}
For $n\in\mathbb N$ and $(r,t)\in J(n)$, with $n-r-t=2k$, 
the number of $\R$-classes in $\DrPBn$ in which each element has $t$ singleton non-transversal upper-kernel classes is equal to
\[
\rho_{nrt} = {n \choose r}{n-r\choose t} (2k-1)!! = \frac{n!}{2^k k! r! t!}.
\]
\end{lemma}

\pf By Corollary \ref{greenPBn}, the number of such $\R$-classes is equal to the number of graphs on vertex set $\bn$ with $r$ vertices of degree $0$, $t$ vertices with a single loop, and the remaining $n-r-t$ vertices of degree $1$.  To specify such a graph, we first choose the vertices of degree $0$ in ${n\choose r}$ ways.  We then choose the vertices with loops in ${n-r\choose t}$ ways.  And finally, we choose the remaining edges in $(n-r-t-1)!!=(2k-1)!!$ ways. \epf

\begin{thm}\label{anrt}
Define a sequence $a_{nrt}$, for $n\in\mathbb N$ and $(r,t)\in J(n)$, by
\begin{align*}
a_{nn0} &= 1 &&\text{for all $n$}\\
a_{n0t} &= a_n &&\text{if $n-t$ is even}\\
a_{nrt} &= a_{n-1,r-1,t} +(n-r-t)a_{n-2,r,t} &&\text{if $n\geq2$ and $1\leq r\leq n-1$,}
\end{align*}
where the sequence $a_n$ is defined in Lemma \ref{lamKn}.  Then for any $n\in\mathbb N$ and $(r,t)\in J(n)$, and with~$\rho_{nrt}$ as in Lemma \ref{rhonrt}:
\ben
\item[\emph{(i)}] $e(\RaPBn)=a_{nrt}$ for any $\al\in\DrPBn$ with $t$ singleton non-transversal upper-kernel classes,
\item[\emph{(ii)}] $e(\DrPBn)=\sum_{t\in I(n-r)}\rho_{nrt}a_{nrt}$,
\item[\emph{(iii)}] $e(\PB_n)=\sum_{(r,t)\in J(n)}\rho_{nrt}a_{nrt}$.
\een
\end{thm}

\pf Again, it suffices to prove (i).  Let $\al\in \DrPBn$.  Re-labelling the points from $\bn$, if necessary, we may assume that
\[
\Lamu(\al) = \textlinegraph{-5.0}{
\redto{10}{11}
\redto78
\dotto{1.4}{2.6}
\dotto{4.4}{5.6}
\dotto{8.4}{9.6}
\draw[red](4,0)arc(270:-270:.2);
\draw[red](6,0)arc(270:-270:.2);
\overt1
\overt3
\overt4
\overt6
\overt7
\overt8
\overt{10}
\overt{11}
\draw(2,0)node[below]{$\underbrace{\phantom{mmmm...}}_r$};
\draw(5,0)node[below]{$\underbrace{\phantom{mmmm...}}_t$};
}\ .
\vspace{-\bigskipamount} 
\]
Let $A_{nrt}$ be the set of all balanced graphs on vertex set $\bn$ with the same red edges as $\Lamu(\al)$.  Again, by Corollary \ref{greenPBn} and Proposition \ref{bijections}, it suffices to show that the numbers $a_{nrt}=|A_{nrt}|$ satisfy the stated recurrence.

It is clear that $a_{nn0}=1$ for all $n$.  If $r=0$ (and $n-t$ is even), then we may complete $\Lamu(\al)$ to a graph from $A_{n0t}$ by adding as many (non-adjacent) blue edges as we like, and adding blue loops to the remaining vertices.  Again, such assignments of blue edges are in one-one correspondence with the involutions of $\bn$, of which there are $a_n$.  Now suppose $n\geq2$ and $1\leq r\leq n-1$.  By inspection of (1--4) in Lemma \ref{Lamalform}, we see that elements of $A_{nrt}$ come in two kinds:
\ben 
\item those for which $1$ is a connected component of its own, and
\item those for which $1$ is an endpoint of an even length alternating path (with no loops).
\een
The proof concludes in similar fashion to that of Theorem \ref{anr}. \epf

\section{Idempotents in diagram algebras}\label{sect:algebras}

Let $\al,\be\in\P_n$.  Recall that the product $\al\be\in\P_n$ is defined in terms of the product graph $\Ga(\al,\be)$.  Specifically, $A$ is a block of $\al\be$ if and only if $A=B\cap(\bn\cup\bn')\not=\emptyset$ for some connected component~$B$ of $\Ga(\al,\be)$.  In general, however, the graph $\Ga(\al,\be)$ may contain some connected components strictly contained in the middle row $\bn''$, and the \emph{partition algebra} $\Pnxi$ is designed to take these components into account.  We write $m(\al,\be)$ for the number of connected components of the product graph $\Ga(\al,\be)$ that are entirely contained in the middle row.  It is important to note (and trivially true) that $m(\al,\be)\leq n$ for all $\al,\be\in\P_n$.
Now let $F$ be a field and fix some $\xi\in F$.  
We denote by $\Pnxi$ the $F$-algebra with basis $\P_n$ and product $\circ$ defined on basis elements $\al,\be\in\P_n$ (and then extended linearly) by
\[
\al\circ\be = \xi^{m(\al,\be)}(\al\be).
\]
If $\al,\be,\ga\in\P_n$, then $m(\al,\be)+m(\al\be,\ga) = m(\al,\be\ga)+m(\be,\ga)$, and it follows that $\Pnxi$ is an associative algebra.  We may also speak of the subalgebras of $\Pnxi$ spanned by $\B_n$ and $\PB_n$; these are the \emph{Brauer} and \emph{partial Brauer algebras} $\Bnxi$ and $\PBnxi$, respectively.  See \cite{HR} for a survey-style treatment of the partition algebras.

In this section, we determine the number of partitions $\al\in\P_n$ such that $\al$ is an idempotent basis element of $\Pnxi$; that is, $\al\circ\al=\al$.  These numbers depend on whether $\xi$ is a root of unity.  As such, we define
\[
M = \begin{cases}
m &\text{if $\xi$ is an $m$th root of unity where $m\leq n$}\\
0 &\text{otherwise.}
\end{cases}
\]
If $\K_n$ is one of $\P_n,\B_n,\PB_n$, we will write
\[
\Exi(\K_n) = \set{\al\in\K_n}{\text{$\al=\al\circ\al$ in $\Knxi$}} \AND \exi(\K_n) =|\Exi(\K_n) |.
\]

\smallskip
\vspace{-\medskipamount} 
\begin{thm}\label{char2}
Let $\al\in\P_n$, and suppose the kernel-classes of $\al$ are $X_1,\ldots,X_k$.  
Then the following are equivalent:
\ben
\item[\emph{(1)}] $\al\in\Exi(\P_n)$,
\item[\emph{(2)}] $\al\in E(\P_n)$ and $\rank(\al)\equiv k\pmod M$,
\item[\emph{(3)}] the following three conditions are satisfied:
\begin{enumerate}
\item[\emph{(i)}] $\al\in\P_{X_1}\oplus\cdots\oplus\P_{X_k}$,
\item[\emph{(ii)}] the restrictions $\al|_{X_i}$ all have rank at most $1$,
\item[\emph{(iii)}] the number of restrictions $\al|_{X_i}$ of rank $0$ is a multiple of $M$.
\end{enumerate}
\een
\end{thm}

\nc{\PXixi}{\P_{X_i}^\xi}

\pf First, note that if $\al\in E(\P_n)$, then Theorem \ref{char} gives $\al=\al_1\oplus\cdots\oplus\al_k$ where $\al_i=\al|_{X_i}\in C(\P_{X_i})$ for each $i$, and $r=\rank(\al)=\rank(\al_1)+\cdots+\rank(\al_k)$ with $\rank(\al_i)\in\{0,1\}$ for each~$i$.  Re-labelling the $X_i$ if necessary, we may suppose that $\rank(\al_1)=\cdots=\rank(\al_r)=1$.  Then the connected components contained entirely in $X''$ in the product graph $\Ga(\al,\al)$ are precisely the sets $X_{r+1}'',\ldots,X_k''$.  So $m(\al,\al)=k-r$.

Now suppose (1) holds.  Then $\al=\al\circ\al=\xi^{m(\al,\al)}(\al^2)$, so $\al=\al^2$ and $m(\al,\al)\in M\mathbb Z$.  Since $\al\in E(\P_n)$, it follows from the first paragraph that $k-r=m(\al,\al)\in M\mathbb Z$, and so (2) holds.

Next, suppose (2) holds.  Since $\al\in E(\P_n)$, Theorem \ref{char} tells us that (i) and (ii) hold.  Write $\al=\al_1\oplus\cdots\oplus\al_k$ where $\al_i=\al|_{X_i}\in C(\P_{X_i})$ for each $i$.  The set $\set{i\in\bk}{\rank(\al_i)=0}$ has cardinality $m(\al,\al)$, which is equal to $k-\rank(\al)$ by the first paragraph.  By assumption, $k-\rank(\al)\in M\mathbb Z$, so (iii) holds.

Finally, suppose (3) holds and write $\al_i=\al|_{X_i}$ for each $i$.  Since $\rank(\al_i)\leq1$ and $X_i$ is a kernel-class of $\al$, it follows that $\al_i\in\P_{X_i}$ is irreducible and so $\al_i\in C(\P_{X_i})$ by Lemma \ref{tech3}.
For each $i\in\bk$, let
\[
l_i = \begin{cases}
0 &\text{if $\rank(\al_i)=1$}\\
1 &\text{if $\rank(\al_i)=0$.}
\end{cases}
\]
Then $l_1+\cdots+l_k$ is a multiple of $M$ by assumption, and $\al_i\circ\al_i=\xi^{l_i}(\al_i^2)=\xi^{l_i}\al_i$ in $\PXixi$ for each~$i$.  But then $\al\circ\al=\xi^{l_1+\cdots+l_k}\al=\al$ so that (1) holds.  \epf

\begin{rem}
If $M=0$, then part (2) of the previous theorem says that $\rank(\al)=k$.  Also, conditions (ii) and (iii) in part (3) may be replaced with the simpler statement that the restrictions $\al|_{X_i}$ all have rank $1$.  If $M=1$, then $\xi=1$ so $\Pnxi$ is the (non-twisted) semigroup algebra of~$\P_n$ and $\Exi(\P_n)=E(\P_n)$; in this case, Theorem \ref{char2} reduces to Theorem \ref{char}.
\end{rem}

We are now ready to give formulae for $\exi(\K_n)$ where $\K_n$ is one of $\P_n,\B_n,\PB_n$.  
It will be convenient to give separate statements depending on whether $M=0$ or $M>0$.  The next result is proved in an almost identical fashion to Theorem \ref{eKn}, relying on Theorem \ref{char2} rather than Theorem \ref{char}.

\bs
\begin{thm}\label{eKnxi=0}
If $M=0$ and $\K_n$ is one of $\P_n,\B_n,\PB_n$, then 
\[
\exi(\K_n) = n! \cdot\sum_{\mu\vdash n} \prod_{i=1}^n \frac{c_1(\K_i)^{\mu_i}}{\mu_i! (i!)^{\mu_i}} .
\]
The numbers $\exi(\K_n)$ satisfy the recurrence:
\[
\exi(\K_0)=1, \quad \exi(\K_n) = \sum_{m=1}^n {n-1\choose m-1} c_1(\K_m)\exi(\K_{n-m}) \quad\text{for $n\geq1$.}
\]
The values of $c_1(\K_n)$ are given in Propositions \ref{cPn}, \ref{cBn} and \ref{cPBn}. \epfres
\end{thm}

Recall that if $\al\in\P_n$ has kernel-classes $X_1,\ldots,X_k$ with $|X_1|\geq\cdots\geq|X_k|$, then the integer partition $\mu(\al)$ is defined to be $(|X_1|,\ldots,|X_k|)$.  For a subset $\Si\sub\P_n$ and an integer partition $\mu\vdash n$, we write
\[
\Eximu(\Si) = \set{\al\in \Exi(\Si)}{\mu(\al)=\mu} \AND \eximu(\Si)=|\Eximu(\al)|.
\]
If $\mu=(m_1,\ldots,m_k)\vdash n$, we call $k$ the \emph{height of $\mu$}, and we write $k=h(\mu)$.  The next result follows quickly from Theorem \ref{char2}.

\bs
\begin{thm}\label{eKnxi>0}
Suppose $M>0$, and let $\K_n$ be one of $\P_n,\B_n,\PB_n$.  Then
\[
\exi(\K_n) = \sum_{\mu\vdash n} \eximu(\K_n).
\]
If $\mu\vdash n$ and $k=h(\mu)$, then
\[
\eximu(\K_n) = \sum_{0\leq r\leq n \atop r\equiv k\;\! (\operatorname{mod} M)} e_\mu(\DrKn).
\]
The values of $e_\mu(\DrKn)$ are given in Theorem \ref{eDrKn}. \epfres
\end{thm}

%
%
%
%

We may also derive recurrences for the values of $\exi(\DrKn)$ in the case $M=0$.  Things get more complicated when $M>0$ since the question of whether or not an element of $E(\K_n)$ belongs additionally to $\Exi(\K_n)$ depends not just on its rank but also on the number of kernel classes.  We will omit the $M>0$ case.

\bs
\begin{thm}\label{eDrKnxi=0}
If $M=0$ and $\K_n$ is one of $\P_n,\B_n,\PB_n$, then the numbers $\exi(\DrKn)$ satisfy the recurrence:
\begin{align*}
\exi(D_n(\K_n)) &= 1 \\ 
\exi(D_0(\K_n)) &= 0 &&\text{if $n\geq1$}\\ 
\exi(\DrKn) &= \sum_{m=1}^n {n-1 \choose m-1} c_1(\K_m)\exi(D_{r-1}(\K_{n-m})) &&\text{if $1\leq r\leq n-1$.}
\end{align*}
The values of $c_1(\K_n)$ are given in Propositions \ref{cPn}, \ref{cBn} and \ref{cPBn}. \epfres
\end{thm}

We now use Theorems \ref{eKnxi=0} and \ref{eKnxi>0} to derive explicit values for $\exi(\B_n)$ and $\exi(\PB_n)$ in the case $M=0$.  In fact, since $c_1(\B_n)=c_1(\PB_n)$ by Propositions \ref{cBn} and \ref{cPBn}, it follows that $\exi(\B_n)=\exi(\PB_n)$ in this case.
These numbers seem to be Sequence A088009 on the OEIS \cite{OEIS}, although it is difficult to understand why.

\bs
\begin{thm}
If $n\in\mathbb N$ and $M=0$, then
\[
\exi(\B_n)=\exi(\PB_n) = \sum_{\mu} \frac{n!}{\mu_1!\mu_3!\cdots\mu_{2k+1}!},
\]
where $k=\lfloor\frac{n-1}2\rfloor$, and the sum is over all integer partitions $\mu=(1^{\mu_1},\ldots,n^{\mu_n})\vdash n$ with $\mu_{2i}=0$ for $i=0,1,\ldots,\lfloor\frac n2\rfloor$. 
The numbers $\exi(\B_n)=\exi(\PB_n)$ satisfy the recurrence:
\[
\epfreseq
\exi(\B_0)=1 \COMMA \exi(\B_n) = \sum_{i=0}^{\lfloor \frac {n-1}2\rfloor} {n-1\choose 2i}(2i+1)! \;\:\!\! \exi(\B_{n-2i-1}) \qquad\text{for $n\geq1$.}
\]
\end{thm}

Theorem \ref{eDrKnxi=0} yields a neat recurrence for the numbers $\exi(\DrBn)=\exi(\DrPBn)$ in the case $M=0$.

\bs
\begin{thm}\label{eDrBnxi=0}
If $M=0$, then the numbers $\exi(\DrBn)=\exi(\DrPBn)$ satisfy the recurrence:
\begin{align*}
\exi(D_n(\B_n)) &= 1 \\ 
\exi(D_0(\B_n)) &= 0 &&\text{if $n\geq1$}\\ 
\epfreseq
\exi(\DrBn) &= \sum_{i=0}^{\lfloor \frac{n-1}2\rfloor} {n-1 \choose 2i} (2i+1)! \;\:\!\! \exi(D_{r-1}(\B_{n-2i-1})) &&\text{if $1\leq r\leq n-1$.}
\end{align*}
\end{thm}

We may also use the methods of Section \ref{sect:diff} to derive recurrences for the number of elements from an $\R$-, $\L$- or $\D$-class of $\B_n$ or $\PB_n$ that satisfy $\al=\al\circ\al$ in $\Pnxi$ with $M=0$.
Recall that $I(n)=\set{r\in\bn^0}{n-r\in2\mathbb Z}$.

\bs
\begin{thm}\label{bnr}
Define a sequence $b_{nr}$, for $n\in\mathbb N$ and $r\in I(n)$, by
\begin{align*}
b_{nn} &= 1 &&\text{for all $n$}\\
b_{n0} &= 0 &&\text{if $n\geq2$ is even}\\
b_{nr} &= b_{n-1,r-1}+(n-r)b_{n-2,r} &&\text{if $1\leq r\leq n-2$.}
\end{align*}
Then for any $n\in\N$ and $r\in I(n)$, and with $M=0$ and $\rho_{nr}$ as in Lemma \ref{rhonr}:
\ben
\item[\emph{(i)}] $\exi(\RaBn)=\exi(\LaBn)=b_{nr}$ for any $\al\in \DrBn$,
\item[\emph{(ii)}] $\exi(\DrBn)=\rho_{nr}b_{nr}$,
\item[\emph{(iii)}] $\exi(\B_n)=\sum_{r\in I(n)}\rho_{nr}b_{nr}$.
\een
\end{thm}

\pf The proof is virtually identical to the proof of Theorem \ref{anr}, except we require $b_{n0}=0$ if $n\geq2$ is even since reduced balanced graphs with $n$ vertices and $n/2$ red (and blue) edges correspond to elements of $D_0(\B_n)$, and these do not belong to $\Exi(\B_n)$ by Theorem \ref{char2}. \epf

\begin{rem}
If $\al\in\PB_n$, then $E(\RaPBn)$ is non-empty if and only if $\RaPBn$ has non-trivial intersection with $\B_n$, in which case $e(\RaPBn)=b_{nr}$ where $r=\rank(\al)$.  (This is because ${\be\in \Exi(\RaPBn)}$ if and only if the connected components of $\Lam(\be)$ are all of type (1) as stated in Lemma~\ref{Lamalform}, in which case $\be\in\Exi(\B_n)$.)  A dual statement can be made concerning $\L$-classes.  It follows that $\exi(\DrPBn)=\exi(\DrBn)=\rho_{nr}b_{nr}$ for all $r\in I(n)$, and $\exi(\PB_n)=\exi(\B_n)=\sum_{r\in I(n)}\rho_{nr}b_{nr}$.
\end{rem}

\section{Calculated values}\label{sect:values}

In this section, we list calculated values of $c_0(\K_n),c_1(\K_n),c(\K_n),e(\K_n),\exi(\K_n)$ where $\K_n$ is one of $\P_n,\B_n,\PB_n$ and where $M=0$.  We also give values of $e(\DrKn)$ and $\exi(\DrKn)$ where $M=0$, and $e(\RaBn)$ and $\exi(\RaBn)$ for $\al\in\DrBn$.

\begin{table}[H]
\begin{center}
\begin{tabular}{|c||rrrrrrrrrrr|}
\hline
$n$ & 0 & $1$ & $2$ & $3$ & $4$ & $5$ & $6$ & $7$ & $8$ & $9$ & $10$ \\
\hline\hline
$c_0(\B_n)$ &   &  0 & 1  & 0  & 6  & 0  &  120  & 0  & 5040  & 0  & 362880 \\
\hline
$c_1(\B_n)$ &   &  1  & 0  & 6  & 0  &  120  & 0  & 5040  & 0  & 362880 & 0 \\
\hline
$c(\B_n)$ &   &  1  & 1  & 6  & 6  &  120  & 120  & 5040  & 5040  & 362880 & 362880 \\
\hline
$e(\B_n)$ &  1 &  1 &  2  & 10  & 40  & 296  &  1936  & 17872  & 164480  & 1820800  & 21442816 \\
\hline
$\exi(\B_n)$ &  1  &  1  &  1  &  7  &  25  &  181  &  1201  &  10291  &  97777  &  1013545  &  12202561  \\
\hline
\end{tabular}
\end{center}
\caption{Calculated values of $c_0(\B_n),c_1(\B_n),c(\B_n),e(\B_n),\exi(\B_n)$ with $M=0$.}
\label{tab:Bn}
\end{table}

\begin{table}[H]
\begin{center}
\begin{tabular}{|c||rrrrrrrrrrr|}
\hline
$n$ & 0 & $1$ & $2$ & $3$ & $4$ & $5$ & $6$ & $7$ & $8$ & $9$ & $10$ \\
\hline\hline
$c_0(\PB_n)$ &   &   1  & 3  & 6  & 30  &  120  & 840  & 5040  & 45360  & 362880 & 3991680 \\
\hline
$c_1(\PB_n)$ &   &  1  & 0  & 6  & 0  &  120  & 0  & 5040  & 0  & 362880 & 0 \\
\hline
$c(\PB_n)$ &   &  2  & 3  & 12  & 30  &  240  & 840  & 10080  & 45360  & 725860 & 3991680 \\
\hline
$e(\PB_n)$ &  1  &  2  &  7  &  38  &  241  &  1922  &  17359  &  180854  &  2092801  &  26851202  &  376371799  \\
\hline
$\exi(\PB_n)$ &  1  &  1  &  1  &  7  &  25  &  181  &  1201  &  10291  &  97777  &  1013545  &  12202561  \\
\hline
\end{tabular}
\end{center}
\caption{Calculated values of $c_0(\PB_n),c_1(\PB_n),c(\PB_n),e(\PB_n),\exi(\PB_n)$ with $M=0$.}
\label{tab:PBn}
\end{table}

\begin{table}[H]
\begin{center}
{\small
\begin{tabular}{|c||rrrrrrrrrrr|}
\hline
$n$ & 0 & $1$ & $2$ & $3$ & $4$ & $5$ & $6$ & $7$ & $8$ & $9$ & $10$ \\
\hline\hline
$c_0(\P_n)$ &   &  1& 3& 15& 119& 1343& 19905& 369113& 8285261& 219627683& 6746244739 \\
\hline
$c_1(\P_n)$ &   &  1  & 5  & 43  & 529  &  8451  & 167397  & 3984807  & 111319257  & 3583777723 & 131082199809 \\
\hline
$c(\P_n)$ &   &  2 &  8 &  58 &  648 &  9794 &  187302 &  4353920 &  119604518  & 3803405406 &  137828444548  \\
\hline
$e(\P_n)$ &  1 &  2& 12& 114& 1512& 25826& 541254& 13479500& 389855014& 12870896154& 478623817564 \\
\hline
$\exi(\P_n)$ &  1  &  1  &  6  &  59  &  807  &  14102  &  301039  &  7618613  & 223586932   &  7482796089  &  281882090283  \\
\hline
\end{tabular}
}
\end{center}
\caption{Calculated values of $c_0(\P_n),c_1(\P_n),c(\P_n),e(\P_n),\exi(\P_n)$ with $M=0$.}
\label{tab:Pn}
\end{table}

\begin{table}[H]
\begin{center}
{
\begin{tabular}{|c||rrrrrrrrrrr|}
\hline
$n\sm r$ & 0 & $1$ & $2$ & $3$ & $4$ & $5$ & $6$ & $7$ & $8$ & $9$ & $10$ \\
\hline\hline
0   &        1             &&&&&&&&&&        \\
1   &              &     1        &&&&&&&&&            \\
2   &        1           &&           1              &&&&&&&&        \\
3   &       &           9           &&           1        &&&&&&&           \\
4   &      9           &&          30           &&           1                  &&&&&&    \\
5   &     &         225           &&          70           &&           1              &&&&&       \\
6   &    225           &&       1575          &&         135           &&           1                  &&&&    \\
7   &   &       11025           &&        6615           &&         231           &&           1                   &&&  \\
8   & 11025           &&      132300           &&       20790           &&         364          &&           1        &&   \\
9   &&      893025          &&      873180          &&       54054          &&         540           &&           1&\\
10\phantom{1}  & 893025           &&    16372125           &&     4054050           &&      122850           &&         765           &&           1   \\
\hline
\end{tabular}
}
\end{center}
\caption{Calculated values of $e(\DrBn)$.}
\label{tab:eDrBn}
\end{table}

\begin{table}[H]
\begin{center}
{
\begin{tabular}{|c||rrrrrrrrrrr|}
\hline
$n\sm r$ & 0 & $1$ & $2$ & $3$ & $4$ & $5$ & $6$ & $7$ & $8$ & $9$ & $10$ \\
\hline\hline
0   &        1             &&&&&&&&&&        \\
1   &              &     1        &&&&&&&&&            \\
2   &        1           &&           1              &&&&&&&&        \\
3   &       &           3           &&           1        &&&&&&&           \\
4   &      3           &&          5           &&           1                  &&&&&&    \\
5   &     &         15           &&          7           &&           1              &&&&&       \\
6   &    15           &&       35          &&         9           &&           1                  &&&&    \\
7   &   &       105           &&        63           &&         11           &&           1                   &&&  \\
8   & 105           &&      315           &&       99           &&         13          &&           1        &&   \\
9   &&      945          &&      693          &&       143          &&         15           &&           1&\\
10\phantom{1}  & 945           &&    3465           &&     1287           &&      195           &&         17           &&           1   \\
\hline
\end{tabular}
}
\end{center}
\caption{Calculated values of $e(\RaBn)=e(\LaBn)$ where $\al\in\DrBn$.}
\label{tab:eLaBn}
\end{table}

\begin{table}[H]
\begin{center}
{
\begin{tabular}{|c||rrrrrrrrrrr|}
\hline
$n\sm r$ & 0 & $1$ & $2$ & $3$ & $4$ & $5$ & $6$ & $7$ & $8$ & $9$ & $10$ \\
\hline\hline
0   &        1             &&&&&&&&&&        \\
1   &              &     1        &&&&&&&&&            \\
2   &        0           &&           1              &&&&&&&&        \\
3   &       &           2           &&           1        &&&&&&&           \\
4   &      0           &&          4           &&           1                  &&&&&&    \\
5   &     &         8           &&          6           &&           1              &&&&&       \\
6   &    0           &&       24          &&         8           &&           1                  &&&&    \\
7   &   &       48           &&        48           &&         10           &&           1                   &&&  \\
8   & 0           &&      192           &&       80          &&         12          &&           1        &&   \\
9   &&      384          &&      480          &&       120          &&         14           &&           1&\\
10\phantom{1}  & 0           &&    1920           &&     168           &&      195           &&         16           &&           1   \\
\hline
\end{tabular}
}
\end{center}
\caption{Calculated values of $\exi(\RaBn)=\exi(\LaBn)=\exi(\RaPBn)=\exi(\LaPBn)$ where $\al\in\DrBn$ and $M=0$.}
\label{tab:exiLaBn}
\end{table}

\begin{table}[H]
\begin{center}
{\small
\begin{tabular}{|c||rrrrrrrrrrr|}
\hline
$n\sm r$ & 0 & $1$ & $2$ & $3$ & $4$ & $5$ & $6$ & $7$ & $8$ & $9$ & $10$ \\
\hline\hline
0   &        1             &&&&&&&&&&        \\
1   &       1      &     1        &&&&&&&&&            \\
2   &        4      &     2      &     1              &&&&&&&&        \\
3   &       16    &      18     &      3      &     1        &&&&&&&           \\
4   &      100    &      88     &     48       &    4       &    1                  &&&&&&    \\
5   &     676     &    860     &    280      &   100       &    5       &    1              &&&&&       \\
6   &    5776    &    6696    &    4020    &     680     &    180       &    6       &    1                 &&&&    \\
7   &   53824   &    76552    &   35196    &   13580     &   1400     &    294      &     7      &     1                   &&&  \\
8   & 583696    &  805568   &   531328  &    131936    &   37240     &   2576    &     448      &     8       &    1        &&   \\
9   & 6864400  &  10765008   &  6159168   &  2571744   &   397656    &   88200    &    4368   &      648      &     9     &      1 &\\
10\phantom{1}  & 90174016 &  141145120 &  101644560  &  32404800   &  9780960   &  1027152  &    187320    &    6960    &     900      &    10           &1    \\
\hline
\end{tabular}
}
\end{center}
\caption{Calculated values of $e(\DrPBn)$.}
\label{tab:eDrPBn}
\end{table}

\begin{table}[H]
\begin{center}
{\tiny
\begin{tabular}{|c||rrrrrrrrrrr|}
\hline
$n\sm r$ & 0 & $1$ & $2$ & $3$ & $4$ & $5$ & $6$ & $7$ & $8$ & $9$ & $10$ \\
\hline\hline
0   &        1             &&&&&&&&&&        \\
1   &       1      &     1        &&&&&&&&&            \\
2   &        4      &     7     &      1              &&&&&&&&        \\
3   &       25     &     70     &     18       &    1        &&&&&&&           \\
4   &      225    &     921     &    331        &  34          & 1                  &&&&&&    \\
5   &     2704   &    15191    &    6880      &   995        &  55        &   1              &&&&&       \\
6   &    41209   &   304442    &  163336     &  29840      &  2345        &  81           &1                 &&&&    \\
7   &   769129   &  7240353   &  4411190    &  958216    &   95760      &  4739         &112         &  1                   &&&  \\
8   & 17139600  & 200542851  & 134522725  &  33395418   &  3992891    &  252770     &   8610       &  148       &    1        &&   \\
9   &447195609  & 6372361738  & 4595689200  & 1267427533 &  174351471   &13274751  & 581196   &14466   &  189 &    1&\\
10\phantom{1}  & 13450200625  & 229454931097   &174564980701  & 52345187560&  8059989925 &  709765413   &37533657 &  1205460  &  2289  &  235  &  1   \\
\hline
\end{tabular}
}
\end{center}
\caption{Calculated values of $e(\DrPn)$.}
\label{tab:eDrPn}
\end{table}

\begin{table}[H]
\begin{center}
{
\begin{tabular}{|c||rrrrrrrrrrr|}
\hline
$n\sm r$ & 0 & $1$ & $2$ & $3$ & $4$ & $5$ & $6$ & $7$ & $8$ & $9$ & $10$ \\
\hline\hline
0   &        1             &&&&&&&&&&        \\
1   &              &     1        &&&&&&&&&            \\
2   &        0           &&           1              &&&&&&&&        \\
3   &       &           6           &&           1        &&&&&&&           \\
4   &     0            &&          24           &&           1                  &&&&&&    \\
5   &     &         120          &&          60           &&           1              &&&&&       \\
6   &    0           &&       1080          &&         120           &&           1                  &&&&    \\
7   &   &       5040           &&        5040           &&         210           &&           1                   &&&  \\
8   &    0        &&      80640           &&       16800           &&         336          &&           1        &&   \\
9   &   &      362880          &&      604800          &&       45360          &&         504           &&           1&\\
10\phantom{1}  &    0        &&    9072000           &&     3024000           &&      105840           &&         720           &&           1   \\
\hline
\end{tabular}
}
\end{center}
\caption{Calculated values of $\exi(\DrBn)=\exi(\DrPBn)$ with $M=0$.}
\label{tab:exiDrBn}
\end{table}

\begin{table}[H]
\begin{center}
{\scriptsize
\begin{tabular}{|c||rrrrrrrrrrr|}
\hline
$n\sm r$ & 0 & $1$ & $2$ & $3$ & $4$ & $5$ & $6$ & $7$ & $8$ & $9$ & $10$ \\
\hline\hline
0   &        1             &&&&&&&&&&        \\
1   &        0      &     1        &&&&&&&&&            \\
2   &        0      &     5     &      1              &&&&&&&&        \\
3   &       0       &     43     &     15      &     1       &&&&&&&           \\
4   &      0        &     529    &     247     &     30    &       1                  &&&&&&    \\
5   &     0         &    8451   &     4795   &      805   &       50   &        1               &&&&&       \\
6   &    0          &   167397  &    108871    &   22710    &    1985     &     75       &    1                 &&&&    \\
7   &   0           &   3984807   &  2855279 &     697501   &    76790  &      4130  &       105      &     1                   &&&  \\
8   & 0& 111319257 &   85458479   & 23520966  &   3070501    &  209930     &   7658     &    140       &    1         &&   \\
9   & 0&3583777723  & 2887069491  & 871103269 &  129732498&   10604811  & 495054 &  13062  & 180  & 1&\\
10\phantom{1}  & 0  & 131082199809  & 109041191431  & 35334384870 &  5843089225&   549314745  & 30842427  & 1046640 &  20910 &  225  &  1   \\
\hline
\end{tabular}
}
\end{center}
\caption{Calculated values of $\exi(\DrPn)$ with $M=0$.}
\label{tab:exiDrPn}
\end{table}

\section*{Acknowledgement}

We thank the editor and anonymous referees for their helpful comments, and particularly for drawing our attention to a number of references.

\footnotesize
\def\bibspacing{-1.1pt}
\bibliography{biblio}
\bibliographystyle{plain}
\end{document}